\magnification=1200
\loadmsam
\loadmsbm
\loadeufm
\loadeusm
\UseAMSsymbols
\input amssym.def

\font\BIGtitle=cmr10 scaled\magstep3
\font\bigtitle=cmr10 scaled\magstep1
\font\boldsectionfont=cmb10 scaled\magstep1
\font\section=cmsy10 scaled\magstep1

\def\scr#1{{\fam\eusmfam\relax#1}}

\def\scrB{{\scr B}}

\def\scrE{{\scr E}}
\def\scrF{{\scr F}}

\def\scrI{{\scr I}}

\def\scrM{{\scr M}}

\def\scrO{{\scr O}}

\def\scrQ{{\scr Q}}
\def\scrS{{\scr S}}

\def\scrV{{\scr V}}

\def\gr#1{{\fam\eufmfam\relax#1}}

	\def\gre{{\gr e}}
	\def\grf{{\gr f}}
	\def\grg{{\gr g}}
	\def\grh{{\gr h}}
	\def\gri{{\gr i}}
	\def\grj{{\gr j}}

	\def\grm{{\gr m}}
	\def\grn{{\gr n}}

\def\db#1{{\fam\msbfam\relax#1}}

 \def\dbF{{\db F}}
\def\dbG{{\db G}}

 \def\dbN{{\db N}}
 
\def\dbQ{{\db Q}}

 \def\dbZ{{\db Z}}

\def\eps{{\varepsilon}}
\def\vph{{\varphi}}

\def\im{\text{Im}}

\def\Ker{\text{Ker}}

\def\Hom{\text{Hom}}
\def\End{\text{End}}
\def\Spec{\text{Spec}}

\def\Lie{\text{Lie}}

\def\leaderfill{\leaders\hbox to 1em
     {\hss.\hss}\hfill}
\def\nspace{\lineskip=1pt\baselineskip=12pt\lineskiplimit=0pt}

\def\finishproclaim{\par\rm
     \ifdim\lastskip<\medskipamount\removelastskip
     \penalty55\medskip\fi}
\def\endproof{$\hfill \square$}
\def\proof{\par\noindent {\it Proof:}\enspace}
\def\references#1{\par
  \centerline{\boldsectionfont References}\medskip
     \parindent=#1pt\nspace}
\def\Ref[#1]{\par\hang\indent\llap{\hbox to\parindent
     {[#1]\hfil\enspace}}\ignorespaces}
\def\Item#1{\par\smallskip\hang\indent\llap{\hbox to\parindent
     {#1\hfill$\,\,$}}\ignorespaces}
\def\ItemItem#1{\par\indent\hangindent2\parindent
     \hbox to \parindent{#1\hfill\enspace}\ignorespaces}

\def\Le{{\mathchoice{\,{\scriptstyle\le}\,}
  {\,{\scriptstyle\le}\,}
  {\,{\scriptscriptstyle\le}\,}{\,{\scriptscriptstyle\le}\,}}}

\def\arrowsim{\,\smash{\mathop{\to}\limits^{\lower1.5pt
  \hbox{$\scriptstyle\sim$}}}\,}

\def\doublemaprights#1#2#3#4{\raise3pt\hbox{$\mathop{\,\,\hbox to     
#1pt{\rightarrowfill}\kern-30pt\lower3.95pt\hbox to
     #2pt{\rightarrowfill}\,\,}\limits_{#3}^{#4}$}}

\def\rightcapdownarrow{\raise9pt\hbox{$\ssize\cap$}\kern-7.75pt
     \Big\downarrow}

\def\rcapmapdown#1{\rightcapdownarrow\kern-1.0pt\vcenter{
     \hbox{$\scriptstyle#1$}}}

\def\rmapdown#1{\Big\downarrow\kern-1.0pt\vcenter{
     \hbox{$\scriptstyle#1$}}}
\def\rightsubsetarrow#1{{\ssize\subset}\kern-4.5pt\lower2.85pt
     \hbox to #1pt{\rightarrowfill}}
\def\longtwoheadedrightarrow#1{\raise2.2pt\hbox to #1pt{\hrulefill}
     \!\!\!\twoheadrightarrow}

\def\Hom{\operatorname{\hbox{Hom}}}

\NoBlackBoxes
\parindent=25pt
\document
\footline={\hfil}
\footline={\hss\tenrm \folio\hss}
\pageno=1
\bigskip

\noindent 
\vskip 1 cm 

\centerline{\BIGtitle Reconstructing p-divisible groups}
\bigskip\smallskip
\centerline{\BIGtitle from their truncations of small level}
\bigskip\medskip
\centerline{\bigtitle Adrian Vasiu, Binghamton University}
\bigskip
\centerline{Final version identical with the galley proofs (modulo style). To appear in Comment. Math. Helv. {\bf 85} (2010).} 
\bigskip
\centerline{{\it Dedicated to the memory of Angela Vasiu}}

\bigskip\noindent
{\bf ABSTRACT.} Let $k$ be an algebraically closed field of characteristic $p>0$. Let $D$ be a $p$-divisible group over $k$. Let $n_D$ be the smallest non-negative integer for which the following statement holds: if $C$ is a $p$-divisible group over $k$ of the same codimension and dimension as $D$ and such that $C[p^{n_D}]$ is isomorphic to $D[p^{n_D}]$, then $C$ is isomorphic to $D$. To the Dieudonn\'e module of $D$ we associate a non-negative integer $\ell_D$ which is a computable upper bound of $n_D$. If $D$ is a product $\prod_{i\in I} D_i$ of isoclinic $p$-divisible groups, we show that $n_D=\ell_D$; if the set $I$ has at least two elements we also show that $n_D\le\max\{1,n_{D_i},n_{D_i}+n_{D_j}-1|i,j\in I,\;j\neq i\}$. We show that we have $n_D\Le 1$ if and only if $\ell_D\Le 1$; this recovers the classification of minimal $p$-divisible groups obtained by Oort. If $D$ is quasi-special, we prove the Traverso truncation conjecture for $D$. If $D$ is $F$-cyclic, we compute explicitly $n_D$. Many results are proved in the general context of latticed $F$-isocrystals with a (certain) group over $k$.



\bigskip\noindent
{\bf Contents}

\medskip
\smallskip\noindent
{1} Introduction \dotfill\ \ \ 1

\smallskip\noindent
{2} Preliminaries \dotfill\ \ \ 7

\smallskip\noindent
{3} The proof of the Main Theorem A \dotfill\ \ \ 9

\smallskip\noindent
{4} Direct applications to p-divisible groups \dotfill\ \ \ 21

\smallskip\noindent
{5} On the Main Theorem B \dotfill\ \ \ 29
\smallskip\noindent
References\dotfill\ \ \ 31

\bigskip
\noindent
{\boldsectionfont 1. Introduction}
\bigskip

Let $p\in\dbN$ be a prime. Let $k$ be an algebraically closed field of characteristic $p$. Let $c,d\in\dbN\cup\{0\}$ be such that $r:=c+d>0$. Let $D$ be a $p$-divisible group over $k$ of codimension $c$ and dimension $d$. The height of $D$ is $r$. Let $n_D\in\dbN\cup\{0\}$ be the smallest number for which the following statement holds: if $C$ is a $p$-divisible group of codimension $c$ and dimension $d$ over $k$ such that $C[p^{n_D}]$ is isomorphic to $D[p^{n_D}]$, then $C$ is isomorphic to $D$. We have $n_D=0$ if and only if $cd=0$. For the existence of $n_D$ we refer to [Ma, Ch. III, Sect. 3], [Tr1, Thm. 3], [Tr2, Thm. 1], [Va1, Cor. 1.3], or [Oo2, Cor. 1.7]. For instance, one has the following gross estimate $n_D\le cd+1$ (cf. [Tr1, Thm. 3]). The classical Dieudonn\'e theory says that the category of $p$-divisible groups over $k$ is antiequivalent to the category of Dieudonn\'e modules over $k$. Thus the existence of $n_D$ gets translated into a suitable problem pertaining to Dieudonn\'e modules and thus to a particular type of latticed $F$-isocrystals over $k$ (see Subsection 1.1 below for precise definitions). 

Traverso's truncation conjecture predicts that $n_D\le\min\{c,d\}$, cf. [Tr3, Sect. 40, Conj. 4]. This surprising and old conjecture is known to hold only in few cases (like for supersingular $p$-divisible groups over $k$; see [NV, Thm. 1.2]). To prove different refinements of this conjecture, one needs to have easy ways to compute and estimate $n_D$. Each estimate of $n_D$ represents progress towards the classification of $p$-divisible groups over $k$; implicitly, it represents progress towards the understanding of the ultimate stratifications defined in [Va1, Subsect. 5.3] and (thus also) of the special fibres of all integral canonical models of Shimura varieties of Hodge type. The goal of the paper is to put forward {\it basic principles} that compute either $n_D$ or some very sharp upper bounds of $n_D$. 

For the sake of generality, a great part of this paper will be worked out in the context of latticed $F$-isocrystals with a (certain) group over $k$. 

\medskip\smallskip\noindent
{\bf 1.1. Latticed $F$-isocrystals.} Let $W(k)$ be the ring of Witt vectors
 with coefficients in $k$. Let $B(k)$ be the field of fractions of $W(k)$. Let $\sigma$ be the Frobenius automorphism of $W(k)$ and $B(k)$ induced from $k$. 

By a {\it latticed $F$-isocrystal} over $k$ we mean a pair $(M,\phi)$, where $M$ is a free $W(k)$-module of finite rank and $\phi:M[{1\over p}]\arrowsim M[{1\over p}]$ is a $\sigma$-linear automorphism. We recall that if $\phi(M)\subseteq M$, then the pair $(M,\phi)$ is called an {\it $F$-crystal} over $k$. We also recall that if $pM\subseteq \phi(M)\subseteq M$, then the pair $(M,\phi)$ is called a {\it Dieudonn\'e module} over $k$ and $\vartheta:=p\phi^{-1}:M\to M$ is called the Verschiebung map of $(M,\phi)$. 

The composite of $W(k)$-linear maps endows $\End(M)$ with a natural structure of a $W(k)$-algebra (and thus also of a Lie algebra over $W(k)$). We denote also by $\phi$ the
$\sigma$-linear automorphism of $\End(M[{1\over p}])$ that takes
$e\in \End(M[{1\over p}])$ to $\phi(e):=\phi\circ
e\circ\phi^{-1}$. Let $G_{B(k)}$ be a connected subgroup of $\pmb{GL}_{M[{1\over p}]}$ such that its Lie algebra $\Lie(G_{B(k)})$ is left invariant by $\phi$ i.e., we have $\phi(\Lie(G_{B(k)}))=\Lie(G_{B(k)})$. Let $G$ be the schematic closure of $G_{B(k)}$ in $\pmb{GL}_M$. The triple $(M,\phi,G)$ is called a {\it latticed $F$-isocrystal with a group} over $k$, cf. [Va1, Def. 1.1 (a)]. Let $\grg:=\Lie(G_{B(k)})\cap\End(M)$; it is a Lie subalgebra of $\End(M)$ which as a $W(k)$-submodule is a direct summand. If $G$ is smooth over $\Spec(W(k))$, then $\grg=\Lie(G)$. Let $n_G\in\dbN\cup\{0\}$ be the {\it $i$-number} of $(M,\phi,G)$ introduced in [Va1, Def. 3.1.4]. Thus $n_G$ is the smallest non-negative integer for  which the following statement holds:

\medskip
$\bullet$ {\it If $g\in G(W(k))$ is congruent to $1_M$ modulo $p^{n_G}$, then there exists $h\in G(W(k))$ which is an isomorphism between $(M,g\phi,G)$ and $(M,\phi,G)$ (equivalently, between $(M,g\phi)$ and $(M,\phi)$). In other words, we have $hg\phi h^{-1}=\phi$ (equivalently, $hg\phi(h)^{-1}=1_M$).}  

\medskip
In [Va1] we developed methods that provide good upper bounds of $n_G$ (see [Va1, Subsubsect. 3.1.3 and Ex. 3.1.5]). The methods used exponential maps and applied to all possible types of affine, integral group schemes $G$ over $\Spec(W(k))$. But when the type of $G$ is simple (like when $G$ is $\pmb{GL}_M$), then one can obtain significantly better bounds. This idea was exploited to some extent in [Va1, Sect. 3.3] and it is brought to full fruition in this paper. Accordingly, in the whole paper we will work under the following assumption:

\medskip\noindent
{\bf 1.1.1. Assumption.} {\it We have $M\neq 0$, the $W(k)$-submodule $\grg$ of $\End(M)$ is a $W(k)$-subalgebra of $\End(M)$ (and not only a Lie subalgebra of $\End(M)$), and (thus) $G$ is the group scheme over $\Spec(W(k))$ of invertible elements of $\grg$.}

\medskip
Typical cases we have in mind: (i) $G$ is either $\pmb{GL}_M$ or a parabolic subgroup scheme of $\pmb{GL}_M$, (ii) $G$ is the centralizer in $\pmb{GL}_M$ of a semisimple $W(k)$-subalgebra of $\End(M)$, and (iii) $\grg$ is $W(k)1_M\oplus\grn$, with $\grn$ a nilpotent subalgebra (without unit) of $\End(M)$.

\medskip\noindent
{\bf 1.1.2. Newton polygon slopes.} Dieudonn\'e's classification of $F$-isocrystals over $k$ (see [Di, Thms. 1 and 2], [Ma, Ch. 2, Sect. 4], [De], etc.) implies that we have a direct sum decomposition $M[{1\over p}]=\bigoplus_{\alpha\in\dbQ} W(\alpha)$
that is left invariant by $\phi$ and that has the property that all Newton polygon slopes of $(W(\alpha),\phi)$ are $\alpha$. We recall that if $m\in\dbN$ is the smallest number such that $m\alpha\in\dbZ$, then there exists a $B(k)$-basis for $W(\alpha)$ which is formed by elements fixed by $p^{-m\alpha}\phi^m$. One says that $(M,\phi)$ is {\it isoclinic} if there exists a rational number $\alpha$ such that we have $M[{1\over p}]=W(\alpha)$. We consider the direct sum decomposition into $B(k)$-vector spaces
$$\End(M[{1\over p}])=L_+\oplus L_0\oplus L_-$$ 
that is left invariant by $\phi$ and such that all Newton polygon slopes of $(L_+,\phi)$ are positive, all Newton polygon slopes of $(L_0,\phi)$ are $0$, and finally all Newton polygon slopes of $(L_-,\phi)$ are negative. We have direct sum decompositions $L_+=\bigoplus_{\alpha,\beta\in\dbQ\,\alpha<\beta} \Hom(W(\alpha),W(\beta))$, $L_-=\bigoplus_{\alpha,\beta\in\dbQ\,\alpha<\beta} \Hom(W(\beta),W(\alpha))$, and $L_0=\bigoplus_{\alpha\in\dbQ} \End(W(\alpha))$. Thus both $L_+$ and $L_-$ are nilpotent subalgebras (without unit) of $\End(M)$. 

We have $L_0=\End(M[{1\over p}])$ if and only if $(M,\phi)$ is isoclinic. 

\medskip\smallskip\noindent
{\bf 1.2. Level modules and torsions.} We define
$$O_+:=\{x\in\End(M)\cap L_+|\phi^q(x)\in\End(M)\cap L_+\;\forall q\in\dbN\}=\cap_{q\in\dbN\cup\{0\}} \phi^{-q}(\End(M)\cap L_+).$$ 
Let $A_0:=\{e\in \End(M)|\phi(e)=e\}$ be the
$\db{Z}_p$-algebra of endomorphisms of $(M,\phi)$. Let $O_0$ be the $W(k)$-span
of $A_0$; it is a $W(k)$-subalgebra of $\End(M)\cap L_0$. We have identities 
$$O_0=A_0\otimes_{\dbZ_p} W(k)=\cap_{q\in\dbN\cup\{0\}} \phi^q(\End(M)\cap L_0)=\cap_{q\in\dbN\cup\{0\}} \phi^{-q}(\End(M)\cap L_0).$$
We also define 
$$O_-:=\{x\in\End(M)\cap L_-|\phi^{-q}(x)\in\End(M)\cap L_-\;\forall q\in\dbN\}=\cap_{q\in\dbN\cup\{0\}} \phi^{q}(\End(M)\cap L_-).$$

As all Newton polygon slopes of $(L_+,\phi)$ are positive, for each $x\in L_+$ the sequence $(\phi^q(x))_{q\in\dbN}$ of elements of $L_+$ converges to $0$ in the $p$-adic topology. This implies that there exists $s\in\dbN$ such that $p^sx\in O_+$. Thus we have $O_+[{1\over p}]=L_+$. As $O_+$ is a $W(k)$-submodule of the finitely generated $W(k)$-module $\End(M)$, we conclude that $O_+$ is a lattice of $L_+$. A similar argument shows that $O_0$ and $O_-$ are lattices of $L_0$ and $L_-$ (respectively). We have the following relations $\phi(O_+)\subseteq O_+$, $\phi(O_0)=O_0=\phi^{-1}(O_0)$, $\phi^{-1}(O_-)\subseteq O_-$, $L_+L_0+L_0L_+\subseteq L_+$, $L_0L_0\subseteq L_0$, and $L_0L_-+L_0L_-\subseteq L_-$. These relations imply that:

\medskip
{\bf (i)} Both $O_+$ and $O_-$ are left and right $O_0$-modules.

\smallskip
{\bf (ii)} The direct sum $O_+\oplus O_0$ (resp. $O_0\oplus O_-$) is a $W(k)$-subalgebra of $\End(M)$ that has $O_+$ (resp. $O_-$) as a nilpotent, two-sided ideal. 
\medskip
Let $O:=O_+\oplus O_0\oplus O_-$; it is a lattice of $\End(M)[{1\over p}]$ contained in $\End(M)$. In general, $O$ is not a $W(k)$-subalgebra of $\End(M)$ (see Example 2.2). Thus we call $O$ the {\it level module} of $(M,\phi)$. 

Let $O_G:=(\grg\cap O_+)\oplus (\grg\cap O_0)\oplus (\grg\cap O_-)$; it is a lattice of $\grg[{1\over p}]$ contained in $\grg$. We refer to $O_G$ as the {\it level module} of $(M,\phi,G)$. We note down that $O=O_{\pmb{GL}_M}$. 

By the {\it level torsion} of $(M,\phi,G)$ we mean the unique number $\ell_G\in\dbN\cup\{0\}$ for which the following inclusions hold 
$$p^{\ell_G}\grg\subseteq O_G\subseteq \grg\leqno (1)$$ 
and which obeys the following two disjoint rules:

\medskip
{\bf (a)} if $\grg=O_G$ and the two-sided ideal of the $W(k)$-algebra $\grg$ generated by $(\grg\cap O_+)\oplus (\grg\cap O_-)$ is not topologically nilpotent, then $\ell_G:=1$;

\smallskip
{\bf (b)} in all other cases, $\ell_G$ is the smallest non-negative integer for which (1) holds.

\medskip\noindent
{\bf 1.2.1. A connection to [Va1].} Let $m_G:=\pmb{T}(\grg,\phi)$ be the Fontaine--Dieudonn\'e torsion of $(\grg,\phi)$ introduced in [Va1, Defs. 2.2.2 (a) and (b)]. We recall that $m_G$ is the smallest non-negative integer with the property that there exists a $W(k)$-submodule $\grm$ of $\grg$ which contains $p^{m_G}\grg$ and for which the pair $(\grm,\phi)$ is a Fontaine--Dieudonn\'e $p$-divisible object over $k$ in the sense of loc. cit. One has a direct sum decomposition $(\grm,\phi)=\bigoplus_{j\in J} (\grm_j,\phi)$ such that each pair $(\grm_j,\phi)$ is an elementary Fontaine--Dieudonn\'e $p$-divisible object over $k$. The pair $(\grm_j,\phi)$ is a special type of isoclinic latticed $F$-isocrystals over $k$ that are definable over $\dbF_p$; let $\alpha_j\in\dbQ$ be the Newton polygon slope of $(\grm_j,\phi)$. One basic property of $(\grm_j,\phi)$ is: if $\alpha_j>0$ (resp. $\alpha_j=0$ or $\alpha_j<0$), then we have $\phi(\grm_j)\subseteq\grm_j$ (resp. $\phi(\grm_j)=\grm_j$ or $\phi^{-1}(\grm_j)\subseteq\grm_j$). Thus if $\alpha_j>0$ (resp. $\alpha_j=0$ or $\alpha_j<0$), then we have $\grm_j\subseteq\grg\cap O_+$ (resp. $\grm_j\subseteq\grg\cap O_0$ or $\grm_j\subseteq\grg\cap O_-$). This implies that $\grm\subseteq O_G$. Therefore we have $\ell_G\le m_G$ except in the case when $\grg=O_G=\grm$ and $\ell_G=1$. This implies that $\ell_G\le\max\{1,m_G\}$. In general, $\ell_G$ can be smaller than $m_G$ (see Example 2.2). 

\medskip\noindent
{\bf 1.2.2. Example.} We assume that all Newton polygon slopes of $(\grg,\phi)$ are $0$. Then we have $O_G=\grg\cap O_0$ and $\ell_G$ is the smallest non-negative integer such that we have inclusions $p^{\ell_G}\grg\subseteq O_G\subseteq\grg$. As the $W(k)$-module $\grg$ is a direct summand of $\End(M)$, we have $O_G=\grg\cap O_0=\grg[{1\over p}]\cap O_0$. This implies that $\phi(O_G)=O_G$ and therefore $O_G$ has a $W(k)$-basis formed by elements of $\grg\cap A_0$. Thus $(O_G,\phi)$ is a  Fontaine--Dieudonn\'e $p$-divisible object over $k$; therefore $\ell_G=m_G$. 

\medskip
Our first main goal is to prove (see Section 3) the following Theorem.

\medskip\smallskip\noindent
{\bf 1.3. Main Theorem A.} {\it We recall that $(M,\phi,G)$ is a latticed $F$-isocrystal with a group over $k$ and that we work under Assumption 1.1.1.

\medskip
{\bf (a)} Then we have an inequality $n_G\Le {\ell_G}$.

\smallskip
{\bf (b)} Assume that $(M,\phi)$ is a direct sum of isoclinic latticed $F$-isocrystals over $k$. Then we have  $n_{\pmb{GL}_M}=\ell_{\pmb{GL}_M}$.}

\medskip
We neither know nor expect examples with $n_{\pmb{GL}_M}<\ell_{\pmb{GL}_M}$. Our second main goal is to apply the Main Theorem A to study $p$-divisible groups over $k$.

\medskip\smallskip\noindent
{\bf 1.4. First applications to $p$-divisible groups.} Let $D$ and $n_D$ be as in the beginning paragraph of the paper. We say that $D$ is isoclinic if its (contravariant) Dieudonn\'e module is isoclinic. If $(M,\phi)$ is the Dieudonn\'e module of $D$, then let
$$\ell_D:=\ell_{\pmb{GL}_M}\in\dbN\cup\{0\}.$$
We call $\ell_D$ the {\it level torsion} of $D$. The following elementary Lemma is our starting point for calculating and estimating $n_D$.

\medskip\noindent
{\bf 1.4.1. Lemma.} {\it We assume that $(M,\phi)$ is the Dieudonn\'e module of $D$. Then we have $n_D=n_{\pmb{GL}_M}$.} 

\medskip
See [Va1, Lem. 3.2.2 and Cor. 3.2.3] and [NV, Thm. 2.2 (a)] for two proofs of Lemma 1.4.1 (the second proof is not stated in the language of latticed $F$-isocrystals with a group). Accordingly, we call $n_D$ the {\it $i$-number} (i.e., the isomorphism number) of $D$. Based on Lemma 1.4.1, we have the following Corollary of the Main Theorem A.

\medskip\noindent
{\bf 1.4.2. Basic Corollary.} {\it  For each non-trivial $p$-divisible group $D$ over $k$ we have $n_D\le\ell_D$. If $D$ is a direct sum of isoclinic $p$-divisible groups over $k$,  then we have $n_D=\ell_D$.} 

\medskip
The inequality $n_D\le\ell_D$ was first checked for the isoclinic case in [Va1, Ex. 3.3.5]. 

\medskip\noindent
{\bf 1.4.3. Proposition.} {\it We assume that $D=\prod_{i\in I} D_i$ is a direct sum of at least two isoclinic $p$-divisible groups over $k$. Then we have the following basic estimate 
$$n_D\le\max\{1,n_{D_i},n_{D_i}+n_{D_j}-1|i,j\in I, j\neq i\}.$$} 
\indent
Proposition 1.4.3 is proved in Subsection 4.5. Example 4.6.2 shows that in general, Proposition 1.4.3 is optimal. The next Proposition (proved in Subsection 4.7) describes the possible range of variation of $n_D$ and $\ell_D$ under isogenies.

\medskip\noindent
{\bf 1.4.4. Proposition.} {\it Let $D\twoheadrightarrow \tilde D$ be an isogeny between non-trivial $p$-divisible groups over $k$. Let $\kappa\in\dbN\cup\{0\}$ be the smallest number such that $p^{\kappa}$ annihilates the kernel of this isogeny. Then we have $n_D\le\ell_D\le\ell_{\tilde D}+2\kappa$. Thus, if $\tilde D$ is a direct sum of isoclinic $p$-divisible groups, then we have $n_D\le\ell_D\le n_{\tilde D}+2\kappa$.}

\medskip
In general, the constant $2\kappa$ of Proposition 1.4.4 is optimal (see Example 4.7.1).

\medskip\smallskip\noindent
{\bf 1.5.  Minimal and quasi-special types.} Let $\scrB=\{e_1,\ldots,e_{r}\}$ be a $W(k)$-basis for $M$. Let $\pi$ be an arbitrary permutation of the set $J_{r}:=\{1,\ldots,r\}$. Let $(M,\phi_{\pi})$ be the Dieudonn\'e module over $k$ with the property that for each $s\in\{1,\ldots,d\}$ we have $\phi_{\pi}(e_s)=pe_{\pi(s)}$ and for each $s\in\{d+1,\ldots,d+c\}$ we have $\phi_{\pi}(e_s)=e_{\pi(s)}$. Let $C_{\pi}$ be a $p$-divisible group over $k$ whose Dieudonn\'e module is $(M,\phi_{\pi})$. For a cycle $\pi_i=(e_{s_1},\ldots,e_{s_{r_i}})$ of $\pi$, let $c_i$ and $d_i=r_i-c_i$ be the number of elements of the sets $\{s_1,\ldots,s_{r_i}\}\cap\{d+1,\ldots,d+c\}$ and $\{s_1,\ldots,s_{r_i}\}\cap\{1,\ldots,d\}$ (respectively), and let $\alpha_i:={{d_i}\over {r_i}}\in\dbQ\cap [0,1]$.

\medskip\noindent
{\bf 1.5.1. Definitions.} We recall that $c$ and $d$ are non-negative integers such that $r:=c+d>0$, that $D$ is a $p$-divisible group over $k$ of codimension $c$ and dimension $d$, and that $J_r=\{1,\ldots,r\}$.

\medskip
{\bf (a)} We say that $D$ is {\it $F$-cyclic} (resp. {\it $F$-circular}), if there exists a permutation $\pi$ (resp. an $r$-cycle permutation $\pi$) of $J_{r}$ such that $D$ is isomorphic to $C_{\pi}$. 

\smallskip
{\bf (b)} We say that $\pi$ is a {\it minimal permutation}, if the following condition holds:

\medskip
{\bf (*)} for each cycle $\pi_i=(e_{s_1},\ldots,e_{s_{r_i}})$ of $\pi$ and for all $q\in\dbN$ and $u\in\{1,\ldots,r_i\}$, we have $\phi_{\pi}^q(e_{s_u})=p^{[q\alpha_i]+\eps_q(s_u)}e_{\pi^q(s_u)}$ for some number $\eps_q(s_u)\in\{0,1\}$. 

\smallskip
{\bf (c)} We say that $D$ is {\it minimal}, if there exists a minimal permutation $\pi$ of $J_{r}$ such that $D$ is isomorphic to $C_{\pi}$.

\smallskip
{\bf (d)} A non-trivial truncated Barsotti--Tate group $B$ of level $1$ over $k$ is called {\it minimal}, if there exists a $p$-divisible group $\tilde D$ over $k$ such that $\tilde D[p]$ is isomorphic to $B$ and $n_{\tilde D}\le 1$.

\smallskip
{\bf (e)} Let $m:=g.c.d.\{c,d\}\in\dbN$ and let $(d_1,r_1):=({d\over m},{r\over m})$. We say that $D$ is {\it isoclinic quasi-special} (resp. {\it isoclinic special}), if we have $\phi^r(M)=p^dM$ (resp. we have $\phi^{r_1}(M)=p^{d_1}M$). We say that $D$ is {\it quasi-special} (resp. {\it special}), if it is a direct sum of isoclinic quasi-special (resp. isoclinic special) $p$-divisible groups over $k$. 

\medskip
The terminology $F$-cyclic and $F$-circular is suggested by Definition 1.2.4 (c) in [Va2]. The terminology minimal $p$-divisible groups and minimal truncated Barsotti--Tate groups of level $1$ is the one used in [Oo3] and [Oo4]. It is easy to check that the above definitions of minimal $p$-divisible groups over $k$ and of minimal truncated Barsotti--Tate groups of level $1$ over $k$ are equivalent to the ones used in [Oo3, Subsect. 1.1] (this also follows from the Main Theorem B below). Moreover $D$ is minimal if and only if $D[p]$ is minimal, cf. Main Theorem B below. The terminology special (see (e)) is as in [Ma, Ch. III, Sect. 2]. If $D$ is $F$-cyclic, then it is also quasi-special but it is not necessarily special (see Lemma 4.2.4 (a) and Example 4.7.1). The class of isomorphism classes of quasi-special $p$-divisible groups of codimension $c$ and dimension $d$ over $k$, is a finite set (see Lemma 4.2.4 (b)); this result recovers and refines slightly [Ma, Ch. III, Sect. 3, Thm. 3.4].

A systematic approach to $C_{\pi}$'s was started in [Va2] and [Va3] using the language of Weyl groups (the role of a permutation $\pi$ of $J_r$ is that one of a representative of the Weyl group of $\pmb{GL}_M$ with respect to its maximal torus that normalizes $W(k)e_s$ for all $s\in J_r$); for instance, we proved that for two permutations $\pi_1,\pi_2$ of $J_{r}$, the $p$-divisible groups $C_{\pi_1}$ and $C_{\pi_2}$ are isomorphic if and only if $C_{\pi_1}[p]$ and $C_{\pi_2}[p]$ are isomorphic (cf. [Va3, Thm. 1.3 (a) and Fact 4.3.1]). The $p$-divisible groups $C_{\pi}$ are also studied in [Oo4] using the language of cyclic words in the variables $\phi$ and $\vartheta$. We note down that in the condition (*), it suffices to consider natural numbers $q$ which are at most equal to the order of $\pi_i$. Thus we view (b) and (d) as a more practical form of [Oo4, Sect. 4]. 

In Subsection 4.6 we prove the following Theorem. 

\medskip\noindent
{\bf 1.5.2. Theorem.} {\it We assume that the non-trivial $p$-divisible group $D$ is quasi-special (for instance, $D$ is $F$-cyclic or special). Then $D$ is a direct sum of isoclinic $p$-divisible groups over $k$ and thus we have $n_D=\ell_D$. Moreover, we have an inequality $n_D\le\min\{c,d\}$ i.e., the Traverso truncation conjecture holds for $D$.} 

\medskip
The proof of the inequality part of Theorem 1.5.2 relies on Proposition 1.4.3 and on an explicit formula for $n_D$ (see property 4.6 (ii); if $D$ is $F$-cyclic, see also Scholium 4.6.1). 

The importance of minimal $p$-divisible groups stems from the following Theorem to be proved in Subsection 5.1.  

\medskip\smallskip\noindent
{\bf 1.6. Main Theorem B.} {\it Let $D$ be a non-trivial $p$-divisible group over $k$. Then the following three statements are equivalent:

\medskip
{\bf (a)} we have $\ell_D\Le 1$;

\smallskip
{\bf (b)} we have $n_D\Le 1$ (equivalently, $D[p]$ is minimal);

\smallskip
{\bf (c)} the $p$-divisible group $D$ over $k$ is minimal.}

\medskip
The implication $(c)\Rightarrow (b)$ was first checked for the isoclinic case in [Va1, Ex. 3.3.6] and for the general case in [Oo3, Thm. 1.2]. A great part of [Oo4] is devoted to the proof of the equivalence between (b) and (c), cf. [Oo4, Thm. B]. 

\bigskip
\noindent
{\boldsectionfont 2. Preliminaries}
\bigskip 

Let $(M,\phi)$ be a latticed $F$-isocrystal over $k$. In this Section we include simple properties that pertain to $(M,\phi)$. Let $M^*:=\Hom(M,W(k))$. 

The notations $p$, $k$, $c$, $d$, $r=c+d$, $D$, $n_D$, $W(k)$, $B(k)$, $(M,\phi,G)$, $\grg$, $M[{1\over p}]=\bigoplus_{\alpha\in\dbQ} W(\alpha)$, $L_+$, $L_0$, $L_-$, $O_+$, $A_0$, $O_0$, $O_-$, $O_G$, $\ell_G$, $\ell_D$, $J_r=\{1,\ldots,r\}$, $(M,\phi_{\pi})$, and $C_{\pi}$ introduced in Section 1, will be used throughout the paper. Let $D^{\text{t}}$ be the $p$-divisible group over $k$ which is the Cartier dual of $D$.  For $m\in\dbN$, let $W_m(k):=W(k)/p^mW(k)$. 

All finitely generated $W(k)$-modules and all finite dimensional $B(k)$-vector spaces are endowed with the $p$-adic topology. As in Subsection 1.2, in the whole paper we keep the following order: first $+$, next $0$, and last $-$.

\medskip\smallskip\noindent
{\bf 2.1. Duals and homs.} Let $\phi:M^*[{1\over p}]\arrowsim M^*[{1\over p}]$ be the $\sigma$-linear automorphism that takes $f\in M^*[{1\over p}]$ to $\sigma\circ f\circ\phi^{-1}\in M^*[{1\over p}]$. The latticed $F$-isocrystal $(M^*,\phi)$ is called the dual of $(M,\phi)$, cf. [Va1, Subsect. 2.1]. The canonical identification $\End(M)=M\otimes_{W(k)} M^*$ defines an identification $(\End(M),\phi)=(M,\phi)\otimes (M^*,\phi)$ of latticed $F$-isocrystals over $k$. If $(M,\phi)$ is the Dieudonn\'e module of $D$, then $(M^*,p\phi)$ is the Dieudonn\'e module of $D^{\text{t}}$. Let $(M_1,\phi_1)$ and $(M_2,\phi_2)$ be two latticed $F$-isocrystals over $k$. Let $\phi_{12}:\Hom(M_1,M_2)[{1\over p}]\arrowsim \Hom(M_1,M_2)[{1\over p}]$ be the $\sigma$-linear automorphism that takes $f\in \Hom(M_1[{1\over p}],M_2[{1\over p}])$ to $\phi_2\circ f\circ\phi_1^{-1}\in \Hom(M_1[{1\over p}],M_2[{1\over p}])$. The latticed $F$-isocrystal $(\Hom(M_1,M_2),\phi_{12})$ over $k$ is called the hom of $(M_1,\phi_1)$ and $(M_2,\phi_2)$. Thus $(M^*,\phi)$ is the hom of $(M,\phi)$ and $(W(k),\sigma)$. The dual of $(\Hom(M_1,M_2),\phi_{12})$ is $(\Hom(M_2,M_1),\phi_{21})$ (here $\phi_{21}$ is defined similarly to $\phi_{12}$). Thus the dual of $(\End(M),\phi)$ is $(\End(M),\phi)$ itself. 

If $\scrB$ is a $W(k)$-basis for $M$, let $\scrB^*:=\{x^*|x\in\scrB\}$ be the dual $W(k)$-basis for $M^*$. Thus for $x,y\in\scrB$, we have $x^*(y)=\delta_{xy}$. For $q\in\dbZ$ and $x,y\in\scrB$, let $a_q(x,y)\in B(k)$ be such that we have  $\phi^q(x)=\sum_{y\in\scrB} a_q(x,y)y$. We have $\phi^q(x^*)=\sum_{y\in\scrB} \sigma^q(a_{-q}(y,x))y^*$ and therefore $\phi^{-q}(x^*)=\sum_{y\in\scrB} \sigma^{-q}(a_q(y,x))y^*$. This implies that:

\medskip
{\bf (*)}  if $s\in\dbZ$, then we have $p^s\phi^q(M)\subseteq M$ (i.e., $p^sa_q(x,y)\in W(k)$ for all $x,y\in\scrB$) if and only if we have $p^s\phi^{-q}(M^*)\subseteq M^*$ (i.e., $p^sa_q(y,x)\in W(k)$ for all $x,y\in\scrB$).  
        
\medskip
The set $\{x\otimes y^*|x,y\in\scrB\}$ is a $W(k)$-basis for $\End(M)=M\otimes_{W(k)} M^*$.

\medskip\smallskip\noindent
{\bf 2.2. Example.} We assume that we have a direct sum decomposition $M=W(k)\oplus N$ such that $\phi$ acts on $W(k)$ as $\sigma$ does, we have $\phi(N)\subseteq N$, and $(N,\phi)$ is isoclinic of Newton polygon slope $\gamma\in (\dbQ\cap (0,\infty))\setminus\dbZ$. We have a direct sum decomposition of latticed $F$-isocrystals over $k$
$$(\End(M),\phi)=(\End(N),\phi)\oplus (N,\phi)\oplus (N^*,\phi)\oplus (W(k),\sigma).$$
The $W(k)$-span of the product $N^*N$ (taken inside $\End(M)$) is $\End(N)$. As $\phi(N)\subseteq N$ and $N^*\subseteq\phi^{-1}(N^*)$, we have $N\subseteq O_+$ and $N^*\subseteq O_-$. As $\gamma\notin\dbZ$, we have $O\cap\End(N)\subsetneqq \End(N)$. Thus $O_+O_-\nsubseteq O$. Therefore $O$ is not a $W(k)$-subalgebra of $\End(M)$.  

We take $G$ such that $\grg$ is the $W(k)$-subalgebra $N\oplus W(k)1_N$ of $\End(M)$. As $N\subseteq O_+$ and $1_M\in O_0$, we have $\grg=O_G=(\grg\cap O_+)\oplus (\grg\cap O_0)$ and $(\grg\cap O_+)$ is a nilpotent, two-sided ideal of the $W(k)$-algebra $\grg$. Thus $\ell_G=0$, cf. the rule 1.2 (a). If the pair $(N,\phi)$ is not a Dieudonn\'e--Fontaine $p$-divisible object over $k$, then the Dieudonn\'e--Fontaine torsion $m_G$ of $(\grg,\phi)$ is positive (and in fact it can be any natural number).

\medskip\smallskip\noindent
{\bf 2.3. Lemma.} {\it We assume that $n_{\pmb{GL}_M}=0$. Then there exists an integer $s$ such that we have $\phi(M)=p^sM$. Thus $\phi(\End(M))=\End(\phi(M))=\End(M)$ and therefore we have $O_0=\End(M)$ and $\ell_{\pmb{GL}_M}=0$.}

\medskip
\proof
Let $q\in\dbN$. By induction on $q$ we show that the Lemma holds if the rank $r$ of $M$ is at most $q$. If $q=1$ and $r=1$, then the Lemma is obvious. The passage from $q$ to $q+1$ goes as follows. We can assume that $r=q+1$. By multiplying $\phi$ with $p^{-s}$ for some $s\in\dbZ$, we can assume that $\phi(M)$ is a $W(k)$-submodule of $M$ that contains a direct summand of $M$ of rank at least $1$. Let $\tilde x\in M\setminus pM$ be such that $\phi(\tilde x)\in M\setminus pM$. Let $g_{\tilde x}\in\pmb{GL}_M(W(k))$ be such that $g_{\tilde x}\phi(\tilde x)=\tilde x$. As $n_{\pmb{GL}_M}=0$, $(M,\phi)$ is isomorphic to $(M,g_{\tilde x}\phi)$. Thus there exists $x\in M\setminus pM$ such that $\phi(x)=x$. Let $M_0$ be the $W(k)$-submodule of $M$ generated by elements fixed by $\phi$; it is a direct summand of $M$ which contains $x$. 

If $M_0=M$, then we are done as $\phi(M)=M$. Thus to end the proof it suffices to show that the assumption that $M_0\neq M$ leads to a contradiction. Let $M_1:=M/M_0$ and let $\phi_1:M_1\to M_1$ be the $\sigma$-linear endomorphism induced by $\phi$. For each element $g_1\in\pmb{GL}_{M_1}(W(k))$ there exists an element $g\in \pmb{GL}_M(W(k))$ that fixes $M_0$ and that maps naturally to $g_1$. As $(M,g\phi)$ and $(M,\phi)$ are isomorphic and due to the definition of $M_0$, we easily get that $(M_1,\phi_1)$ and $(M_1,g_1\phi_1)$ are isomorphic. Thus the $i$-number of $(M_1,\phi_1,\pmb{GL}_{M_1})$ is $0$. As the rank of $M_1$ is less than $q+1$, by induction we get that there exists a natural number $s_1$ such that $\phi(M_1)=p^{s_1}M_1$ (we have $s_1\neq 0$, due to the definition of $M_0$). Let $z_1\in M_1\setminus pM_1$ be such that $\phi_1(z_1)=p^{s_1}z_1$. Let $\tilde z\in M$ be such that it maps naturally to $z_1$. We have $\phi(\tilde z)-p^{s_1}\tilde z\in M_0$. Let $\tilde y\in M_0$ be such that $\phi(\tilde y)-p^{s_1}\tilde y=-\phi(\tilde z)+p^{s_1}\tilde z$. If $z:=\tilde z+\tilde y$, then we have $\phi(z)=p^{s_1}z$. As $z$ maps naturally to $z_1\in M_1\setminus pM_1$, the $W(k)$-module $M_0\oplus W(k)z$ is a direct summand of $M$. Let $g_{xz}\in\pmb{GL}_M(W(k))$ be such that it permutes $x$ and $z$, it normalizes $M_0\oplus W(k)z$, and it acts identically on $(M_0\oplus W(k)z)/(W(k)x\oplus W(k)z)$ and on $M/(M_0\oplus W(k) z)$. The Newton polygon slopes of $(M, g_{xz}\phi)$ are $0$, ${s_1\over 2}$, and $s_1$. As the Newton polygon slopes of $(M,\phi)$ are $0$ and $s_1$ and as $s_1\in\dbN$, we get that $(M,\phi)$ and $(M,g_{xz}\phi)$ are not isomorphic. This contradicts the equality $n_{\pmb{GL}_M}=0$.\endproof

\medskip\smallskip\noindent
{\bf 2.4. Lemma.} {\it Let $x\in\End(M)$ be such that for all $q\in\dbN$ (resp. for all $q\in -\dbN$) we have $\phi^q(x)\in\End(M)$. Then we have $x\in O_+\oplus O_0$ (resp. we have $x\in O_0\oplus O_-$).}

\medskip
\proof
We will prove only the non-negative part of the Lemma as the non-positive part of it is proved in the same way. Thus we assume that we have $\phi^q(x)\in\End(M)$ for all $q\in\dbN$. We write $x=x_++x_0+x_-$, where $x_+\in L_+$, $x_0\in L_0$, and $x_-\in L_-$. There exists a number $s\in\dbN$ such that $p^sx_+\in O_+$ and $p^sx_0\in O_0$. Thus $\phi^q(p^sx_+)\in O_+\subseteq\End(M)$ and $\phi^q(p^sx_0)\in O_0\subseteq\End(M)$. We easily get that we have $p^s\phi^q(x_-)\in\End(M)$ for all $q\in\dbN$. This implies that $x_-=0$ (as all Newton polygon slopes of $(L_-,\phi)$ are negative). Thus $x=x_++x_0$. The sequence $(\phi^q(x_+))_{q\in\dbN}$  converges to $0$ (as all Newton polygon slopes of $(L_+,\phi)$ are positive).  Thus there exists $\tilde q\in\dbN$ such that $y_+:=\phi^{\tilde q}(x_+)\in O_+$. Let $y:=\phi^{\tilde q}(x)$ and $y_0:=y-y_+=\phi^{\tilde q}(x_0)\in\End(M)\cap L_0$. As for each $q\in\dbN$ we have $\phi^q(y_+)\in O_+\subseteq\End(M)$ and $\phi^q(y)\in\End(M)$, we also have $\phi^q(y_0)\in\End(M)$. Thus $y_0\in O_0$. Therefore $x_0=\phi^{-\tilde q}(y_0)\in \phi^{-\tilde q}(O_0)=O_0$. This implies that for all $q\in\dbN\cup\{0\}$ we have $\phi^q(x_0)\in O_0\subseteq\End(M)$. Thus for all $q\in\dbN\cup\{0\}$ we have $\phi^q(x_+)=\phi^q(x)-\phi^q(x_0)\in\End(M)$ i.e., $x_+\in O_+$. Therefore $x=x_++x_0\in O_+\oplus O_0$.\endproof

\medskip\smallskip\noindent
{\bf 2.5. Invertible elements.} In this Subsection we recall basic properties of invertible elements of $W(k)$-subalgebras of $\End(M)$. Let $\grh$ be a $W(k)$-subalgebra of $\End(M)$.  

\medskip
{\bf (i)} If $x\in \grh$ has an inverse $x^{-1}$ in $\End(M)$ (i.e., if $x\in \grh\cap\pmb{GL}_M(W(k))$), then the determinant of $x$ is an invertible element of $W(k)$ and therefore from the Cayley--Hamilton theorem we get that $x^{-1}$ is a polynomial in $x$ with coefficients in $W(k)$; thus $x^{-1}\in\grh$ (i.e., $x$ is an invertible element of $\grh$). 

\smallskip
{\bf (ii)} Each invertible element of $\grh$ is also an invertible element of any other $W(k)$-subalgebra of $\End(M)$ that contains  $\grh$.

\smallskip
{\bf (iii)} If we have a direct sum decomposition $\grh=\grn\oplus\grh_0$ such that $\grh_0$ is a $W(k)$-subalgebra of $\grh$ and $\grn$ is a nilpotent, two-sided ideal of the $W(k)$-algebra $\grh$, then we have a short exact sequence $1\to 1_M+\grn\to \grh\cap\pmb{GL}_M(W(k))\to \grh_0\cap\pmb{GL}_M(W(k))\to 1$ which splits and which is defined by the rule: if $x\in\grn$ and $y\in\grh_0$ are such that $x+y\in \grh\cap\pmb{GL}_M(W(k))$, then the image of $x+y$ in $\grh_0\cap\pmb{GL}_M(W(k))$ is $y$. 

\smallskip
{\bf (iv)} We recall that a two-sided ideal $\gri$ of the $W(k)$-algebra $\grh$ is called {\it topologically nilpotent} if for all $m\in\dbN$ there exists $\tilde m\in \dbN$ such that we have an inclusion $\gri^{\tilde m}\subseteq p^m\grh$ (this implies that $\cap_{m\in\dbN} \gri^m=0$). If $x\in\gri$, then the element $1_M+\sum_{m=1}^{\infty} (-x)^m\in 1_M+\gri$ is well defined and is the inverse of $1_M+x$. This implies that an element of $\grh$ is invertible if and only if its image in $\grh/\gri$ is an invertible element of $\grh/\gri$.

\bigskip
\noindent
{\boldsectionfont 3. The proof of the Main Theorem A}
\bigskip

In this Section we prove the Main Theorem A (see Subsections 3.4 and 3.5). We begin by introducing certain $W(k)$-algebras and group schemes over $\Spec(W(k))$ and by presenting basic properties of them  (see Subsections 3.1 and 3.2). In Subsection 3.3 we list simple properties of isomorphism classes of certain latticed $F$-isocrystals over $k$. All these properties play a key role in Subsection 3.4. In Subsections 3.6 and 3.7 we include two remarks as well as a more general variant of Theorem 1.3 (b).

\medskip\smallskip\noindent
{\bf 3.1. Group schemes of invertible elements.} Let $\grh_+:=\grg\cap O_+$ and $\grh_0:=\grg\cap O_0$. Let $\grh:=(\grh_+\oplus \grh_0)+ p^{\ell_G}\grg$. As $O_0$ and $O_+\oplus O_0$ are $W(k)$-algebras and as $p^{\ell_G}\grg$ is a two-sided ideal of the $W(k)$-algebra $\grg$, it is easy to see that $\grh$ is a $W(k)$-subalgebra of $\grg$. Let $\grh_-:=\grh\cap O_-$. As $p^{\ell_G}\grg\subseteq O_G$ (see (1)), we have a direct sum decomposition
$$\grh=(\grh_+\oplus \grh_0)+ p^{\ell_G}\grg=\grh_+\oplus\grh_0\oplus\grh_-.\leqno (2)$$
Let $\Pi_+:\grh\to\grh$ be the projection on $\grh_+$ along $\grh_0\oplus\grh_-$. 

\indent
Let $H$, $H_{+0}$, and $H_{0}$ be the affine group scheme over $\Spec(W(k))$ of invertible elements of $\grh$, $\grh_+\oplus\grh_0$, and $\grh_0$ (respectively). Due to Subsection 2.5 (ii), we have a sequence 
$$H_0(W(k))\leqslant H_{+0}(W(k))\leqslant H(W(k))\leqslant G(W(k))\leqno (3a)$$
of subgroups. As $\grh_+$ and $\grh_-$ are nilpotent subalgebras (without unit) of $\grh$, we have
$$1_M+\grh_+\leqslant H_{+0}(W(k))\;\;\text{and}\;\; 1_M+\grh_-\leqslant H(W(k)).\leqno (3b)$$
From Subsection 2.5 (iii) we get that we have a natural split short exact sequence
$$1\to 1_M+\grh_+\to H_{+0}(W(k))\to H_0(W(k))\to 1.\leqno (3c)$$
Based on (2), for each element $h\in H(W(k))$ we can write uniquely 
$$h=1_M+a(h)+b(h)+c(h),$$
where $a(h)\in\grh_+$, $b(h)\in\grh_0$, and $c(h)\in\grh_-$. We have $a(h)=\Pi_+(h)$. 

\medskip\noindent
{\bf 3.1.1. The ideal $\gri$.} If $\ell_G=0$, let $\gri$ be the two-sided ideal of the $W(k)$-algebra $\grh=\grg$ generated by $\grh_+$ and $\grh_-$. If $\ell_G\ge 1$, let $\gri:=\grh_++\grh_-+p^{\ell_G}\grg$. We check that $\gri$ is a topologically nilpotent, two-sided  ideal of $\grh$. If $\ell_G=0$, this is so by the very definitions (see rules 1.2 (a) and (b)). We assume that $\ell_G\ge 1$. This implies that $p^{\ell_G}\grg$ is a topologically nilpotent, two-sided ideal of $\grh$. As $\grh_+$ is a nilpotent, two-sided ideal of $\grh_+\oplus\grh_0$, its image in $\grh/p^{\ell_G}\grg=\grh_+\oplus\grh_0/[(\grh_+\oplus\grh_0)\cap p^{\ell_G}\grg]$ is a nilpotent, two-sided ideal. Thus $p^{\ell_G}\grg+\grh_+$ is a topologically nilpotent, two-sided ideal of $\grh$. As $\grh_-$ is a nilpotent, two-sided ideal of $\grh_0\oplus\grh_-$ and as $\grh_0\oplus\grh_-$ surjects onto $\grh/(p^{\ell_G}\grg+\grh_+)=\grh_0/[\grh_0\cap (p^{\ell_G}\grg+\grh_+)]$, the image of $\grh_-$ in $\grh/(p^{\ell_G}\grg+\grh_+)$ is a nilpotent, two-sided ideal. From the last two sentences, we get that $\gri=p^{\ell_G}\grg+\grh_++\grh_-$ is a topologically  nilpotent, two-sided  ideal of $\grh$.

\medskip\noindent
{\bf 3.1.2. Fact.} {\it For each element $h=1_M+a(h)+b(h)+c(h)\in H(W(k))$, we have $1_M+b(h)\in H_0(W(k))$. Therefore also $1_M+a(h)+b(h)\in H_{0+}(W(k))$.}

\medskip
\proof 
As $1_M+b(h)$ and $h$ are congruent modulo $\gri$, the first part of the Fact follows from Subsection 2.5 (iv). The last part of the Fact follows from $(3c)$.\endproof

\medskip\noindent
{\bf 3.1.3. On $H_0$.} As $\grh_0=\grg\cap O_0$, we have $\phi(\grh_0)=\grh_0$ (see Example 1.2.2). Let $\grh_{0\dbZ_p}$ be the $\dbZ_p$-subalgebra of $\grh_0$ formed by elements fixed by $\phi$. Let $H_{0\dbZ_p}$ be the affine group scheme over $\Spec(\dbZ_p)$ of invertible elements of $\grh_{0\dbZ_p}$. The group scheme $H_{0\dbZ_p}$ is a $\dbZ_p$-structure of $H_0$ and thus the Frobenius automorphism $\sigma$ acts naturally on $H_0(W(k))=H_{0\dbZ_p}(W(k))$: for $*\in H_0(W(k))$ we have $\sigma(*)=\phi(*)$. The scheme $H_{0\dbZ_p}$  is an open subscheme of the vector group scheme over $\Spec(\dbZ_p)$ defined by $\grh_{0\dbZ_p}$ (viewed only as a $\dbZ_p$-module). Thus the affine, smooth group scheme $H_{0\dbZ_p}$ has connected fibres. 

\medskip\noindent
{\bf 3.1.4. Lemma.} {\it Let $\grf_+$ be a $W(k)$-submodule of $\grh_+$. Let $\grj:=\grf_+\oplus\grh_0\oplus\grh_-$. We consider the following three conditions: 

\medskip
{\bf (i)} we have $\Pi_+(\grf_+\grh_-+\grh_-\grf_+)\subseteq\grf_+$;

\smallskip
{\bf (ii)} the $W(k)$-module $\grf_+$ is a left and right $\grh_0$-module;

\smallskip
{\bf (iii)} we have $\grf_+^2\subseteq\grf_+$ (i.e., $\grf_+$ is an algebra).

\medskip\noindent
Then the following three properties hold:

\medskip
{\bf (a)} Conditions (i) and (ii) hold if and only if $\grj$ is a left and right $\grh_0\oplus\grh_-$-module.

\smallskip
{\bf (b)} Conditions (ii) and (iii) hold if and only if $\grf_+\oplus\grh_0$ is a $W(k)$-subalgebra of $\grh_+\oplus\grh_0$.

\smallskip
{\bf (c)} The three conditions (i) to (iii) hold if and only if $\grj$ is a $W(k)$-subalgebra of $\grh$.}

\medskip
\proof
As $\grh_0\oplus\grh_-$ is a $W(k)$-subalgebra of $\grh$, $\grj$ is a left and right $\grh_0\oplus\grh_-$-module if and only if we have $\grf_+\grh_0+\grh_0\grf_++\grf_+\grh_-+\grh_-\grf_+\subseteq\grj$. We have $\grf_+\grh_-+\grh_-\grf_+\subseteq\grj$ if and only if (i) holds. As $\grh_+$ is a left and right $\grh_0$-module, we have $\grf_+\grh_0+\grh_0\grf_+\subseteq\grj$ if and only if $\grf_+\grh_0+\grh_0\grf_+\subseteq\grf_+$ and thus if and only if (ii) holds. Part (a) follows from the last three sentences. As $\grh_+$ is an algebra and a left and right $\grh_0$-module, we have $(\grf_+\oplus\grh_0)^2\subseteq\grf_+\oplus\grh_0$ if and only if $\grf_+^2+\grf_+\grh_0+\grh_0\grf_+\subseteq\grf_+$ and thus if and only if conditions (ii) and (iii) hold. Thus (b) holds. Part (c) follows from (a) and (b).\endproof

\medskip\smallskip\noindent
{\bf 3.2. Subalgebras.} In this Subsection we list several subalgebras of $\grh$.

\medskip\noindent
{\bf 3.2.1. Frobenius filtration of $\grh_+$.} For $i\in\dbN\cup\{0\}$ let 
$$\grh_{+,i}:=\grh_+\cap\phi^i(\grg)=\grh_+\cap\phi^i(\grg\cap L_+).$$
We have $\grh_{+,0}=\grh_+$, $\phi(\grh_{+,i})\subseteq\grh_{+,i+1}\subseteq\grh_{+,i}$, and each $\grh_{+,i}$ is a $W(k)$-module and a nilpotent algebra. As $\phi^i(\grh_0)=\grh_0$ (see Subsubsection 3.1.3) and as $\grh_+$ is a left and right $\grh_0$-module, $\grh_{+,i}$ is also a  left and right $\grh_0$-module. As all Newton polygon slopes of $(\grh_+,\phi)$ are positive, we have $\cap_{i=0}^{\infty} \grh_{+,i}=0$. Thus $(\grh_{+,i})_{i\in\dbN\cup\{0\}}$ is a decreasing, separated, and exhaustive filtration of $\grh_+$ to be called the Frobenius filtration.

\medskip\noindent
{\bf 3.2.2. The Theta operations.} We assume that $\grh_+\neq 0$. Let $\scrM(\grh_+)$ be the set of $W(k)$-submodules of $\grh_+$ endowed with the pre-order relation defined by inclusions. We consider the increasing operators $\Theta,\Theta_{\text{a}},\Theta_{\text{s}}:\scrM(\grh_+)\to\scrM(\grh_+)$ that take $\grf_+\in\scrM(\grh_+)$ to 
$$\Theta(\grf_+):=\grf_+^2+\Pi_+(\grf_+\grh_-+\grh_-\grf_+)+\phi(\grf_+)\in\scrM(\grh_+),$$
$\Theta_{\text{a}}(\grf_+):=\grf_+^2+\Pi_+(\grf_+\grh_-+\grh_-\grf_+)$, and $\Theta_{\text{s}}(\grf_+):=\Pi_+(\grf_+\grh_-+\grh_-\grf_+)$. We have identities $\Theta(\grf_+)=\Theta_{\text{a}}(\grf_+)+\phi(\grf_+)$ and $\Theta_{\text{a}}(\grf_+)=\grf_+^2+\Theta_{\text{s}}(\grf_+)$. The lower right indices $\text{a}$ and $\text{s}$ stand for algebraic and slope module (respectively), as suggested by Lemma 3.1.4 (a) and (c). For $i\in\dbN\cup\{0\}$ let 
$$\grf_i:=\Theta^i(\grh_+).$$ 
As $\Theta$ is increasing and as $\Theta(\grh_+)\subseteq\grh_+$, we have $\grf_{i+1}=\Theta(\grf_i)\subseteq\grf_i\subseteq \grh_+$.

\medskip\noindent
{\bf 3.2.3. Lemma.} {\it We assume that $\grh_+\neq 0$. Let $i\in\dbN\cup\{0\}$.  Then $\gre_i:=\grf_i\oplus\grh_0\oplus\grh_-$ is a $W(k)$-subalgebra of $\grh$.}

\medskip
\proof
We use induction on $i$. For $i=0$ we have $\gre_0=\grh$ and thus the Lemma holds. The passage from $i$ to $i+1$ goes as follows. We check that the three conditions (i) to (iii) of Lemma 3.1.4 hold for $\grf_+:=\grf_{i+1}$. As $\grf_{i+1}\subseteq\grf_i$, we have $\Theta_{\text{s}}(\grf_{i+1})\subseteq\Theta_{\text{s}}(\grf_i)\subseteq\Theta(\grf_i)=\grf_{i+1}$. Thus condition 3.1.4 (i) holds. To check that condition 3.1.4 (ii) holds, it suffices to show that each one of the following four elements $\grf_i^2$, $\Pi_+(\grf_i\grh_-)$, $\Pi_+(\grh_-\grf_i)$, and $\phi(\grf_i)$ of $\scrM(\grh_+)$ are left and right $\grh_0$-modules; we will only check that they are left $\grh_0$-modules as the arguments for checking that they are right $\grh_0$-modules are entirely the same. As $\grf_i$ is a left $\grh_0$-module, $\grf_i^2$ is also a left $\grh_0$-module. We have $\grh_0\Pi_+(\grf_i\grh_-)=\Pi_+(\grh_0\grf_i\grh_-)=\Pi_+(\grf_i\grh_-)$ (the last equality as $\grf_i$ is a left $\grh_0$-module). We have $\grh_0\Pi_+(\grh_-\grf_i)=\Pi_+(\grh_0\grh_-\grf_i)=\Pi_+(\grh_-\grf_i)$ (the last equality as $\grh_-$ is a left $\grh_0$-module). We have $\grh_0\phi(\grf_i)=\phi(\grh_0)\phi(\grf_i)=\phi(\grh_0\grf_i)=\phi(\grf_i)$. Thus condition 3.1.4 (ii) holds. As $\grf_{i+1}^2\subseteq\grf_i^2\subseteq\grf_{i+1}$, condition 3.1.4 (iii) also holds. Thus $\gre_{i+1}$ is a $W(k)$-subalgebra of $\grh$, cf. Lemma 3.1.4 (c). This ends the induction.\endproof

\medskip\noindent
{\bf 3.2.4. Lemma.} {\it We assume that $\grh_+\neq 0$. Let $i\in\dbN\cup\{0\}$. Then $\Theta_{\text{s}}(\grh_{+,i})\subseteq \grh_{+,i}$. Thus $\grh_{+,i}\oplus\grh_0\oplus\grh_-$ is a $W(k)$-subalgebra of $\grh$.} 

\medskip 
\proof
Let $x\in \grh_{+,i}$ and $y\in\grh_-$. As $z:=-xy+\Pi_+(xy)\in\grh_0\oplus\grh_-$, we have $\tilde z:=\phi^{-i}(z)\in \grg\cap (O_0\oplus O_-)$. As $x\in\grh_{+,i}$ and $y\in\grh_-\subseteq \grg\cap O_-$, we have $\tilde x:=\phi^{-i}(x)\in\grg\cap L_+$ and $\tilde y:=\phi^{-i}(y)\in\grg\cap L_-$. Thus $\Pi_+(xy)=z+xy=\phi^i(\tilde z+\tilde x\tilde y)\in\phi^i(\grg)$ i.e., $\Pi_+(xy)\in\grh_+\cap\phi^i(\grg)=\grh_{+,i}$. A similar argument shows that $\Pi_+(yx)\in\grh_{+,i}$. Thus $\Theta_{\text{s}}(\grh_{+,i})\subseteq \grh_{+,i}$. As $\grh_{+,i}$ is an algebra and a left and right $\grh_0$-module (see Subsubsection 3.2.1), from Lemma 3.1.4 (c) we get that $\grh_{+,i}\oplus\grh_0\oplus\grh_-$ is a $W(k)$-subalgebra of $\grh$.\endproof

\medskip\noindent
{\bf 3.2.5. Lemma.} {\it We assume that $\grh_+\neq 0$. Then $\grf_{\text{a},\infty}:=\cap_{i=0}^{\infty} \Theta_{\text{a}}^i(\grh_+)$ is $0$.}

\medskip
\proof
Let $\gri$ be as in Subsubsection 3.1.1. Let $\grn_0$ be the topologically nilpotent, two-sided ideal of the $W(k)$-algebra $\grh_0$ such that we have $\gri=\grh_+\oplus\grn_0\oplus\grh_-$. We will check by induction on $q\in\dbN$ that $\grf_{\text{a},\infty}\subseteq\gri^q+\grn_0+\grh_-$. As $\grf_{\text{a},\infty}\subseteq\grh_+\subseteq\gri$, the basis of the induction holds. The passage from $q$ to $q+1$ goes as follows. Let $i\in\dbN$ be such that $\Theta_{\text{a}}^i(\grh_+)\subseteq \gri^q+\grn_0+\grh_-\subseteq\gri$. We have 
$$\Theta_{\text{a}}^{i+1}(\grh_+)=(\Theta_{\text{a}}^i(\grh_+))^2+\Theta_{\text{s}}(\Theta_{\text{a}}^i(\grh_+))\subseteq (\gri^q+\grn_0+\grh_-)^2+\Theta_{\text{s}}(\grh_+\cap (\gri^q+\grn_0+\grh_-))$$
$$\subseteq \gri^{q+1}+\grn_0+\grh_-+\Theta_{\text{s}}(\grh_+\cap (\gri^q+\grn_0+\grh_-)).$$ Let $x\in\grh_+\cap (\gri^q+\grn_0+\grh_-)\subseteq\gri$ and $y\in\grh_-\subseteq\gri$. We have $\Pi_+(xy)-xy\in\grn_0\oplus\grh_-$ and $xy\in\gri^{q+1}+(\grn_0+\grh_-)\grh_-\subseteq\gri^{q+1}+\grn_0+\grh_-$. Thus $\Pi_+(xy)=[\Pi_+(xy)-xy]+xy\in \gri^{q+1}+\grn_0+\grh_-$. A similar argument shows that $\Pi_+(yx)\in \gri^{q+1}+\grn_0+\grh_-$. From the last two sentences we get that $\Theta_{\text{s}}(\grh_+\cap (\gri^q+\grn_0+\grh_-))\subseteq \gri^{q+1}+\grn_0+\grh_-$. We conclude that $\Theta_{\text{a}}^{i+1}(\grh_+)\subseteq \gri^{q+1}+\grn_0+\grh_-$. This implies that $\grf_{\text{a},\infty}\subseteq\gri^{q+1}+\grn_0+\grh_-$. This ends the induction. 

As $\gri$ is topologically nilpotent, we have $\cap_{q\in\dbN} (\gri^q+\grn_0+\grh_-)\subseteq \cap_{q\in\dbN} (p^q\grh+\grn_0+\grh_-)=\grn_0+\grh_-$. This implies that $\grf_{\text{a},\infty}\subseteq\grn_0+\grh_-$. Thus $\grf_{\text{a},\infty}\subseteq\grh_+\cap (\grn_0+\grh_-)=0$ i.e., $\grf_{\text{a},\infty}=0$.\endproof

\medskip\noindent
{\bf 3.2.6. Lemma.} {\it We assume that $\grh_+\neq 0$. Then $\grf_{\infty}:=\cap_{i=0}^{\infty} \grf_i$ is $0$.}

\medskip
\proof
We show that the assumption that $\grf_{\infty}\neq 0$ leads to a contradiction. As we have inclusions $0\subsetneq \grf_{\infty}\subseteq \grf_i\subseteq \grh_+=\grh_{+,0}$ and as $\cap_{i=0}^{\infty} \grh_{+,i}=0$, there exists a greatest number $i_0\in\dbN\cup\{0\}$ for which there exists $i\in\dbN$ such that we have inclusions $\grf_{\infty}\subseteq\grf_i\subseteq\grh_{+,i_0}$. 

As $\Theta_{\text{s}}(\grh_{+,i_0+1})\subseteq\grh_{+,i_0+1}$ (cf. Lemma 3.2.4) and as $\grh_{+,i_0+1}^2+\phi(\grh_{+,i_0+1})\subseteq\grh_{+,i_0+1}$ (cf. Subsubsection 3.2.1), we have $\Theta(\grh_{+,i_0+1})\subseteq\grh_{+,i_0+1}$. Based on this and the inclusion $\phi(\grh_{+,i_0})\subseteq\grh_{+,i_0+1}$, an easy induction on $j\in\dbN$ shows that the images of $\grf_{i+j}=\Theta^j(\grf_i)$ and $\Theta_{\text{a}}^j(\grf_i)$ in $\grh_{+,i_0}/\grh_{+,i_0+1}$ coincide. Let $j_0\in\dbN$ be such that we have $\Theta_{\text{a}}^{j_0}(\grh_+)\subseteq\grh_{+,i_0+1}$, cf. Lemma 3.2.5. Thus the image of $\grf_{i+j_0}$ in $\grh_{+,i_0}/\grh_{+,i_0+1}$ is $0$. Therefore $\grf_{\infty}\subseteq \grf_{i+j_0}\subseteq\grh_{+,i_0+1}$ and this contradicts the choice of $i_0$. Thus $\grf_{\infty}=0$.\endproof 
 
\medskip\smallskip\noindent
{\bf 3.3. Isomorphism properties.} In this Subsection we list properties of the isomorphism classes of those latticed $F$-isocrystals with a group over $k$ which are of the form $(M,g\phi,G)$ with $g\in G(W(k))$. We recall that $\sigma$ acts on $H_0(W(k))$ as $\phi$ does, cf. Subsubsection 3.1.3.

\medskip\noindent
{\bf 3.3.1. Lemma.}  {\it {\bf (a)} We have $H_0(W(k))=\{*^{-1}\phi(*)|*\in H_0(W(k))\}$. 
\medskip
{\bf (b)} If $m\in\dbN$, then $\Ker(H_0(W(k))\to H_0(W_m(k)))=\{*^{-1}\phi(*)|*\in \Ker(H_0(W(k))\to H_0(W_m(k)))\}$.

\smallskip
{\bf (c)} For each $*\in H_{+0}(W(k))$, we have $\phi(*)\in H_{+0}(W(k))$.

\smallskip
{\bf  (d)} Let $\grf_+$ and $\grf_-$ be two left and right $\grh_0$-modules contained in $\grh_+$ and $\grh_-$ (respectively). Let $g\in H(W(k))$ be such that $a(g)\in \grf_+$ and $c(g)\in \grf_-$. Then there exists an element $h_0\in H_0(W(k))$ such that for $g_0:=h_0 g\phi(h_0)^{-1}\in H(W(k))$ we have $a(g_0)\in \grf_+$, $b(g_0)=0$, and $c(g_0)\in \grf_-$.}

\medskip
\proof
As $H_{0\dbZ_p}$ is an affine, smooth group scheme over $\Spec(\dbZ_p)$ whose special fibre is connected (see Subsubsection 3.1.3), (a) and (b) are only the Witt vectors version of Lang theorem for affine, connected, smooth groups over $\dbF_p$; see [NV, Prop. 2.1] and its proof for details. As $\phi(\grh_+)\subseteq\grh_+$ and $\phi(H_0(W(k)))=H_0(W(k))$, from $(3c)$ we get that for each $*\in H_{+0}(W(k))$ we have $\phi(*)\in H_{+0}(W(k))$. Thus (c) holds.

We prove (d). We have $1_M+b(g)\in H_0(W(k))$, cf. Fact 3.1.2. Let $h_0\in H_0(W(k))$ be such that $1_M+b(g)=h_0^{-1}\phi(h_0)$, cf. (a). We have $g_0=h_0a(g)\phi(h_0)^{-1}+h_0[1_M+b(g)]\phi(h_0)^{-1}+h_0c(g)\phi(h_0)^{-1}=1_M+h_0a(g)\phi(h_0)^{-1}+h_0c(g)\phi(h_0)^{-1}$. As $\grf_+$ and $\grf_-$ are left and right $\grh_0$-modules, we have $h_0a(g)\phi(h_0)^{-1}\in \grf_+$ and $h_0c(g)\phi(h_0)^{-1}\in \grf_-$. Therefore $a(g_0)=h_0a(g)\phi(h_0)^{-1}\in \grf_+$ and $c(g_0)=h_0c(g)\phi(h_0)^{-1}\in \grf_-$. Thus (d) holds.\endproof

\medskip\noindent
{\bf 3.3.2. Lemma.} {\it Let $g=1_M+c(g)\in 1_M+\grh_-$. Then there exists an element $h\in G(W(k))\cap (1_M+\grh_-[{1\over p}])$ such that we have $hg\phi(h)^{-1}=1_M$.}

\medskip
\proof 
For $i\in\dbN\cup\{0\}$ let $g_i:=\phi^{-i}(g)=1_M+\phi^{-i}(c(g))$. As $\phi^{-i}(O_-)\subseteq O_-$, we have $\phi^{-i}(\grh_-)\subseteq\grg\cap O_-\subseteq\grg\cap \grh_-[{1\over p}]$ and thus $\phi^{-i}(c(g))$ is a nilpotent element of $\grg$. This implies that $g_{i}=1_M+\phi^{-i}(c(g))$ is an invertible element of $\grg$ i.e., we have $g_{i}\in G(W(k))$. We have $g=g_0$. For $i\in\dbN$ we have $\phi(g_i)=g_{i-1}$. As all Newton polygon slopes of $(L_-,\phi)$ are negative, the sequence $(\phi^{-i}(c(g)))_{i\in\dbN}$ of elements of $\grg\cap O_-$ converges to $0$. This implies that the element $h:=\lim_{i\to\infty} g_ig_{i-1}\cdots g_1\in G(W(k))$
is well defined. We compute that
$$hg\phi(h)^{-1}=\lim_{i\to\infty} g_ig_{i-1}\cdots g_1g \phi(g_1)^{-1}\cdots \phi(g_{i-1})^{-1}\phi(g_i)^{-1}=\lim_{i\to\infty} g_i\cdots g_0g_0^{-1}\cdots g_{i-1}^{-1}$$
is equal to $\lim_{i\to\infty} g_i=1_M$.\endproof

\medskip\smallskip\noindent
{\bf 3.4. Proof of 1.3 (a).} We prove Theorem 1.3 (a). Let $\tilde g\in G(W(k))$ be congruent to $1_M$ modulo $p^{\ell_G}$. As $\tilde g-1_M\in p^{\ell_G}\grg\subseteq\grh$, we have $\tilde g\in\grh$. As $\tilde g\in\pmb{GL}_M(W(k))$, we have $\tilde g\in H(W(k))$ (cf. Subsection 2.5 (i)). Thus to prove Theorem 1.3 (a), it suffices to prove the following stronger statement: 

\medskip
{\bf (*)} {\it for each element $g$ in $H(W(k))$, there exists an element $h_g$ in $G(W(k))$ such that $h_gg\phi(h_g)^{-1}=1_M$}. 

\medskip\noindent
We will first prove the following Lemma. 

\medskip\noindent
{\bf 3.4.1. Lemma.} {\it Let $g\in H(W(k))$. Then there exists an element $h_{+}\in H_{+0}(W(k))$ such that $g_+:=h_{+}g\phi(h_{+})^{-1}\in H(W(k))$ has the property that $a(g_+)=0$.} 

\medskip
\proof
We can assume that $\grh_+\neq 0$. For $i\in \dbN\cup\{0\}$, let $\gre_i=\grf_i\oplus\grh_0\oplus\grh_-$ be the $W(k)$-subalgebra of $\grh$ constructed in Lemma 3.2.3. By induction on $i\in \dbN\cup\{0\}$ we show that there exists $h_i\in H_{+0}(W(k))$ such that $g_i:=h_ig\phi(h_i)^{-1}\in H(W(k))$ has the property that $a(g_i)\in\gre_i$. Taking $h_0=1_M$, we have $g_0=g\in\grh=\gre_0$. Thus the basis of the induction holds. The passage from $i$ to $i+1$ goes as follows. 

We will take $h_{i+1}$ to be a product of the form $h_{i,+}h_{i,0}h_i$. Let $\delta_i\in\dbN\cup\{0\}$ be the greatest number such that we have $b(g_i)\in p^{\delta_i}\grh_0$. Let $h_{i,0}\in\Ker(H_0(W(k))\to H_0(W_{\delta_i}(k)))$ be such that we have $h_{i,0}(1_M+b(g_i))\phi(h_{i,0})^{-1}=1_M$, cf. Lemma 3.3.1 (b). The element $g_{i+1,0}:=h_{i,0}g_i\phi(h_{i,0})^{-1}\in H(W(k))$ has the properties that $a(g_{i+1,0})\in\grf_i$ and $b(g_{i+1,0})=0$, cf. proof of Lemma 3.3.1 (d) applied with $(\grf_+,\grf_-)=(\grf_i,\grh_-)$. Let $h_{i,+}:=1_M-a(g_{i+1,0})\in 1_M+\grh_+\leqslant H_{+0}(W(k))$. We compute that
$$g_{i+1}=h_{i+1}g\phi(h_{i+1})^{-1}=h_{i,+}h_{i,0}g_i\phi(h_{i,0})^{-1}\phi(h_{i,+})^{-1}=h_{i,+}g_{i+1,0}\phi(h_{i,+})^{-1}$$
$$=[1_M-a(g_{i+1,0})][1_M+a(g_{i+1,0})+c(g_{i+1,0})][1_M-\phi(a(g_{i+1,0}))]^{-1}$$
$$=[1_M-a(g_{i+1,0})^2-a(g_{i+1,0})c(g_{i+1,0})+c(g_{i+1,0})][1_M-\phi(a(g_{i+1,0}))]^{-1}.$$ 
As $a(g_{i+1,0})\in\grf_i$, the three elements $-a(g_{i+1,0})^2$, $\Pi_+(-a(g_{i+1,0})c(g_{i+1,0}))$, and $\phi(a(g_{i+1,0}))$ belong to $\Theta(\grf_i)=\grf_{i+1}$. As $\Pi_+(-a(g_{i+1,0})c(g_{i+1,0}))\in\grf_{i+1}$, we get that $-a(g_{i+1,0})c(g_{i+1,0})\in\gre_{i+1}$. As $\gre_{i+1}$ is a $W(k)$-algebra, we conclude that both $1_M-\phi(a(g_{i+1,0}))$ and $[1_M-a(g_{i+1,0})][1_M+a(g_{i+1,0})+c(g_{i+1,0})]$ belong to $\gre_{i+1}$. From Subsection 2.5 (i) we get that $[1_M-\phi(a(g_{i+1,0}))]^{-1}\in\gre_{i+1}$. Thus we have $g_{i+1}\in\gre_{i+1}$. This ends the induction.

Due to Lemma 3.2.6, the sequences $(a(g_i))_{i\in\dbN\cup\{0\}}$ and $(a(g_{i+1,0}))_{i\in\dbN\cup\{0\}}$ of elements of $\grh_+$ converge to $0$. We have $b(g_{i+1})=[1_M-a(g_{i+1,0})][1_M+a(g_{i+1,0})+c(g_{i+1,0})][1_M-\phi(a(g_{i+1,0}))]^{-1}-1_M-a(g_{i+1})-c(g_{i+1})\in\grh_0$. From the last two sentences we easily get that the sequence $(b(g_{i+1}))_{i\in\dbN\cup\{0\}}$ of elements of $\grh_0$ converges to $0$. Thus the sequence $(\delta_i)_{i\in\dbN\cup\{0\}}$ of non-negative integers converges to $\infty$. This implies that the sequence $(h_{i,0})_{i\in\dbN\cup\{0\}}$ of elements of $H_0(W(k))$ converges to $1_M$. As $h_{i,+}=1_M+a(g_{i+1,0})$, the sequence $(h_{i,+})_{i\in\dbN\cup\{0\}}$ converges to $1_M$. Thus the sequence $(h_{i,+}h_{i,0})_{i\in\dbN\cup\{0\}}$ of elements of $H_{+0}(W(k))$ converges also to $1_M$. As $h_{i+1}=h_{i,+}h_{i,0}h_i$, we get that the sequence $(h_i)_{i\in\dbN\cup\{0\}}$ of elements of $H_{+0}(W(k))$ converges to an element $h_+\in H_{+0}(W(k))$. We have $g_+=h_+g\phi(h_+)^{-1}=\lim_{i\to\infty} h_ig\phi(h_i)^{-1}=\lim_{i\to\infty} g_i\in \cap_{i=0}^{\infty} \gre_i=\cap_{i=0}^{\infty} \grf_i\oplus\grh_0\oplus\grh_-$. Thus $g_+\in\grh_0\oplus\grh_-$, cf. Lemma 3.2.6. Therefore $a(g_+)=0$.\endproof

\medskip\noindent
{\bf 3.4.2. End of the proof of 1.3 (a).} Let $g\in H(W(k))$. Let $h_+\in H_{+0}(W(k))$ and $g_+\in H(W(k))$ be as in Lemma 3.4.1. Let $h_0\in H_0(W(k))$ be such that for $g_0:=h_{0}g_+\phi(h_0)^{-1}\in H(W(k))$ we have $a(g_0)=b(g_0)=0$, cf. Lemma 3.3.1 (d) applied with $(\grf_+,\grf_-)=(0,\grh_-)$. Let $h_-\in G(W(k))$ be such that we have $g_0=h_-^{-1}\phi(h_-)$, cf. Lemma 3.3.2. Due to $(3a)$, the element $h_g:=h_-h_0h_+$ belongs to $G(W(k))$. We have $h_gg\phi(h_g)^{-1}=h_-h_0h_+g\phi(h_+)^{-1}\phi(h_0)^{-1}\phi(h_-)^{-1}=h_-h_0g_+\phi(h_0)^{-1}\phi(h_-)^{-1}=h_-g_0\phi(h_-)^{-1}=1_M$. Thus the statement 3.4 (*) holds. This ends the proof of Theorem 1.3 (a).\endproof 

\medskip\noindent
{\bf 3.4.3. Remarks.} {\bf (a)} The proof of Theorem 1.3 (a) can be also worked out using $p^{\ell_G}\grg+\grh_0+\grh_-$ instead of $\grh=p^{\ell_G}\grg+\grh_++\grh_0$.  

\smallskip
{\bf (b)} If $g\in H_{+0}(W(k))$, then $h_g=h_0h_+\in H_{+0}(W(k))$. Thus we have an identity $H_{+0}(W(k))=\{*^{-1}\phi(*)|*\in H_{+0}(W(k))\}$ (to be compared with Lemma 3.3.1 (a)).

\medskip\smallskip\noindent
{\bf 3.5. Proof of 1.3 (b).} We prove Theorem 1.3 (b). We consider a direct sum decomposition 
$$M=\bigoplus_{i\in I}M_i\leqno (4a)$$ 
with the property that for all elements $i$ of the finite set $I$ we have $\phi(M_i[{1\over p}])=M_i[{1\over p}]$ and $(M_i,\phi)$ is isoclinic. For instance, we can take $I$ to be the set of Newton polygon slopes of $(M,\phi)$ and then as each $M_i$ we can take $M\cap W(i)$ (see Subsubsection 1.1.2 for $W(\alpha)$ with $\alpha\in\dbQ$). For each $i\in I$, let $\alpha_i\in\dbQ$ be the unique Newton polygon slope of $(M_i,\phi)$. In Subsubsection 3.5.1 we do not assume that the association $i\to\alpha_i$ is one-to-one.

\medskip\noindent
{\bf 3.5.1. Scholium.} One computes $\ell_{\pmb{GL}_M}$ as follows. For $i\in I$, let $\scrB_i$ be a $W(k)$-basis for $M_i$. Let $\scrB:=\cup_{i\in I} \scrB_i$; it is a $W(k)$-basis for $M$. Let $\scrB^*:=\{x^*|x\in\scrB\}$ be the $W(k)$-basis for $M^*$ which is the dual of $\scrB$ (see Subsection 2.1).

 Due to $(4a)$, we have direct sum decompositions 
$$\End(M)\cap L_+=\bigoplus_{i,j\in I,\alpha_i<\alpha_j} \Hom(M_i,M_j),$$ 
$$\End(M)\cap L_0=\bigoplus_{i\in I} \End(M_i),$$ 
and 
$$\End(M)\cap L_-=\bigoplus_{i,j\in I,\,\alpha_i<\alpha_j} \Hom(M_j,M_i).$$
Thus $\End(M)=(\End(M)\cap L_+)\oplus (\End(M)\cap L_0)\oplus (\End(M)\cap L_-)$ and therefore
$$\End(M)/O=[(\End(M)\cap L_+)/O_+]\oplus [(\End(M)\cap L_0)/O_0]\oplus [(\End(M)\cap L_-)/O_-].\leqno (4b)$$ 
\indent
For $i,j\in I$, $x\in\scrB_i$, and $y\in\scrB_j$, we define a number $\ell(x,y)\in\dbN\cup\{0\}$ via the following two rules: 

\medskip
$\bullet$ if $\alpha_i\ge\alpha_j$, let $\ell(x,y)\in\dbN\cup\{0\}$ be the smallest number such that we have $p^{\ell(x,y)}\phi^q(x\otimes y^*)\in\Hom(M_j,M_i)$ for all $q\in\dbN$;

\smallskip
$\bullet$ if $\alpha_i<\alpha_j$, let $\ell(x,y)\in\dbN\cup\{0\}$ be the smallest number such that we have $p^{\ell(x,y)}\phi^{-q}(x\otimes y^*)\in\Hom(M_j,M_i)$ for all $q\in\dbN$. 

\medskip
Let $\ell_+$, $\ell_0$, $\ell_-\in\dbN\cup\{0\}$ be the smallest numbers such that $p^{\ell_+}$ annihilates $(\End(M)\cap L_+)/O_+$, $p^{\ell_0}$ annihilates $(\End(M)\cap L_0)/O_0$, and $p^{\ell_-}$ annihilates $(\End(M)\cap L_-)/O_-$. As $O_+=\cap_{q\in\dbN\cup\{0\}} \phi^{-q}(\End(M)\cap L_+)=\cap_{q\in\dbN\cup\{0\}} \End(\phi^{-q}(M)\cap L_+)$, $\ell_+$ is the smallest non-negative integer with the property that we have $p^{\ell_+}(\End(M)\cap L_+)\subseteq \phi^{-q}(M)\cap L_+$ for all $q\in\dbN$ (i.e., we have $p^{\ell_+}\phi^q(\End(M)\cap L_+)\subseteq\End(M)$ for all $q\in\dbN$). As $\{x\otimes y^*|x\in\scrB_i,y\in\scrB_j,i,j\in I,\alpha_i>\alpha_j\}$ is a $W(k)$-basis for $\End(M)\cap L_+$, we get that $\ell_+$ is the smallest non-negative integer such that we have $p^{\ell_+}\phi^q(x\otimes y^*)\in\Hom(M_j,M_i)$ for all $q\in\dbN$, all $i,j\in I$ with $\alpha_i>\alpha_j$, and all $x\in\scrB_i$ and $y\in\scrB_j$. Therefore
$$\ell_+:=\max\{\ell(x,y)|x\in\scrB_i,y\in\scrB_j,i,j\in I,\alpha_i>\alpha_j\}.\leqno (5a)$$ 
\noindent
Similar arguments show that 
$$\ell_0=\max\{\ell(x,y)|x,y\in\scrB_i,i\in I\}\leqno (5b)$$
and that
$$\ell_-:=\max\{\ell(x,y)|x\in\scrB_i,y\in\scrB_j,i,j\in I,\alpha_i<\alpha_j\}.\leqno (5c)$$
From $(4b)$ and the very definitions of $\ell_+$, $\ell_0$, and $\ell_-$ we get that $\max\{\ell_+,\ell_0,\ell_-\}\in\dbN\cup\{0\}$ is the smallest number such that $p^{\max\{\ell_+,\ell_0,\ell_-\}}$ annihilates $\End(M)/O$.

Next we define a number $\eps_{\pmb{GL}_M}\in\{0,1\}$ via the following rules. If $O=\End(M)$, let $\eps_{\pmb{GL}_M}:=\ell_{\pmb{GL}_M}$ (cf. rules 1.2 (a) and (b)); we have $\ell_+=\ell_0=\ell_-=0$ and thus $\ell_{\pmb{GL}_M}=\max\{\eps_{\pmb{GL}_M},\ell_+,\ell_0,\ell_-\}$. If $O\neq \End(M)$, let $\eps_{\pmb{GL}_M}:=0$; we have $\ell_{\pmb{GL}_M}=\max\{\ell_+,\ell_0,\ell_-\}$ (cf. rule 1.2 (b)). From the last two sentences and the formulas $(5a)$, $(5b)$, and $(5c)$ we get that, regardless of what $O$ is, we have 
$$\ell_{\pmb{GL}_M}=\max\{\eps_{\pmb{GL}_M},\ell_+,\ell_0,\ell_-\}=\max\{\eps_{\pmb{GL}_M},\ell(x,y)|x,y\in\scrB\} .\leqno (6a)$$
The latticed $F$-isocrystals $(\Hom(M_j,M_i),\phi)$ and $(\Hom(M_i,M_j),\phi)$ are dual to each other (cf. Subsection 2.1) and the dual of the $W(k)$-basis $\{x\otimes y^*|x\in\scrB_i,y\in\scrB_j\}$ of $\Hom(M_j,M_i)$ is the $W(k)$-basis $\{y\otimes x^*|x\in\scrB_i,y\in\scrB_j\}$ of $\Hom(M_i,M_j)$. Based on this, from the property 2.1 (*) we get that for all $i,j\in I$ we have an equality
$$\max\{\ell(x,y)|x\in\scrB_i,y\in\scrB_j\}=\max\{\ell(y,x)|x\in\scrB_i,y\in\scrB_j\}.\leqno (6b)$$ 
\noindent
{\bf 3.5.2. Reduction steps and notations.} Let $\ell:=\ell_{\pmb{GL}_M}$. Based on Theorem 1.3 (a), we have $n_{\pmb{GL}_M}\le\ell$. Thus to prove that $n_{\pmb{GL}_M}=\ell$, it suffices to show that $n_{\pmb{GL}_M}> \ell-1$. If $n_{\pmb{GL}_M}=0$, then $\ell=0$ (see Lemma 2.3) and therefore $n_{\pmb{GL}_M}>\ell-1$. Thus to prove that $n_{\pmb{GL}_M}=\ell$, it suffices to show that  for $\ell\ge 2$ we have $n_{\pmb{GL}_M}> \ell-1$. To check this we can assume that the map $I\to\dbQ$ that takes $l\in I$ to $\alpha_l\in\dbQ$ is injective (i.e., for each element $l\in I$ we have $M_l=M\cap W(\alpha_l)$). 

Let $q$ be the smallest positive integer for which the following two properties hold:
\medskip
{\bf (i)} there exists a $W(k)$-basis $\scrB=\cup_{l\in I} \scrB_l$ for $M$ which is contained in $\cup_{l\in I} M_l$ and for which there exist elements $i,j\in I$, $x\in\scrB_i\subseteq M_i$, and $y\in\scrB_j\subseteq M_j$ such that (cf. $(6a)$ and $(6b)$) we have $\ell(x,y)=\ell$ and $\alpha_j\le\alpha_i$;
\smallskip
{\bf (ii)} we have $e_{q,x,y}:=p^{\ell(x,y)}\phi^q(x\otimes y^*)=p^{\ell}\phi^q(x\otimes y^*)\in\Hom(M_j,M_i)\setminus p\Hom(M_j,M_i)$. 
\medskip\noindent
The existence of $q$ follows from $(6a)$, $(6b)$, and the very definition of the numbers $\ell(x,y)$. 

For $z\in M$, let $a_{z,q}$ be the unique integer such that we have $\phi^q(z)\in p^{a_{z,q}}M\setminus p^{a_{z,q}+1}M$. We can choose the $W(k)$-basis $\scrB=\cup_{l\in I} \scrB_l$ such that we have a direct sum decomposition $M=\bigoplus_{z\in\scrB} W(k)p^{-a_{z,q}}\phi^q(z)$ i.e., we have $\phi^{-q}(M)=\bigoplus_{z\in\scrB} W(k)p^{-a_{z,q}}z$. Let 
$$a_{i,q}:=\min\{a_{z,q}|z\in\scrB_i\}\;\; \text{and}\;\;b_{j,q}:=\max\{a_{z,q}|z\in\scrB_j\}.$$ 
Thus $a_{i,q}$ is the greatest integer such that we have $\phi^q(M_i)\subseteq p^{a_{i,q}}M_i$ and $b_{j,q}$ is the smallest integer such that we have $p^{b_{j,q}}M_j\subseteq \phi^q(M_j)$. The smallest number $s\in\dbN\cup\{0\}$ with the property that $p^s\phi^q(\Hom(M_j,M_i))=p^s\Hom(\phi^q(M_j),\phi^q(M_i))$ is contained in $\Hom(M_j,M_i)$, is equal to $\max\{0,b_{j,q}-a_{i,q}\}$; as $e_{q,x,y}\in \Hom(M_j,M_i)\setminus p\Hom(M_j,M_i)$, we have $s\ge \ell(x,y)=\ell\ge 2$. As $s\le\max\{\ell_+,\ell_0\}\le\ell$, we conclude that $2\le\ell=s=b_{j,q}-a_{i,q}$. It is easy to see that we have $\max\{\ell_+,\ell_0\}\ge a_{y,q}-a_{x,q}\ge \ell(x,y)$, cf. property (ii) for the second inequality. From the last two sentences we get that $a_{x,q}=a_{i,q}$ and $a_{y,q}=b_{j,q}$. Thus we have $\ell=\ell(x,y)=a_{y,q}-a_{x,q}=b_{j,q}-a_{i,q}$. As $\ell=a_{y,q}-a_{x,q}>0$, we have $x\neq y$. 

\medskip\noindent
{\bf 3.5.3. The set $\Lambda$.} Let $\Lambda:=\{w\in M_i\setminus pM_i|a_{w,q}=a_{x,q}\}$; it is the set of those elements $w\in M_i$ for which $p^{-a_{x,q}}\phi^q(w)$ is a direct summand of $M$. Obviously the set $\Lambda$ is stable under multiplication by invertible elements of $W(k)$. For $w\in\Lambda$ let 
$$g_w:=1_M+p^{\ell-1}w\otimes y^*\in\End(M);$$ it is the endomorphism of $M$ that fixes each element $z\in\scrB\setminus\{y\}$ and that takes $y$ to $y+p^{\ell-1}w$. As $\ell\ge 2$, we have $g_w\in\pmb{GL}_M(W(k))$. As each $g_w$ is congruent to $1_M$ modulo  $p^{\ell-1}$, to prove that  $n_{\pmb{GL}_M}>\ell-1$ it suffices to show that there exists an element $w\in\Lambda$ such that the latticed $F$-isocrystals $(M,g_w\phi)$ and $(M,\phi)$ are not isomorphic. We show that the assumption that this is not true leads to a contradiction. This assumption implies that for each element $w\in\Lambda$ there exists an element $h_w\in \pmb{GL}_M(W(k))$ which is an isomorphism between $(M,g_w\phi)$ and $(M,\phi)$. Thus we have $h_wg_w\phi h_w^{-1}=\phi$ i.e., we have 
$$h_wg_w=\phi(h_w).\leqno (7a)$$
We write $h_w=1_M+u_w$, where $u_w\in \End(M)$. Substituting the expressions of $h_w$ and $g_w$ in $(7a)$, we come across the following identity
$$u_w+p^{\ell-1}w\otimes y^*+p^{\ell-1}u_w(w\otimes y^*)=u_w+p^{\ell-1}[w+u_w(w)]\otimes y^*=\phi(u_w)\leqno (7b)$$ 
(here $u_w(w\otimes y^*)$ is the product inside $\End(M)$ of $u_w$ and $w\otimes y^*$). In other words, if $v_w:=w+u_w(w)$ then the pair $(u_w,v_w)$ is a solution of the following equation 
$$U+p^{\ell-1}V\otimes y^*=\Phi(U)\leqno (7c)$$
in variables $U$ and $V$ that can take values in $\End(M)[{1\over p}]$ and $M$ (respectively).

\medskip\noindent
{\bf 3.5.4. Fact.} {\it There exists an isomorphism between $(M,g_w\phi)$ and $(M,\phi)$ defined by an element $\tilde h_w$ of $\pmb{GL}_M(W(k))$ which has the following two properties:

\medskip
{\bf (i)} it acts identically on each $M_l$ with $l\in I\setminus\{i,j\}$ and leaves invariant $M_i$; 

\smallskip
{\bf (ii)} if $i\neq j$, then it acts identically on $M_i$, leaves invariant $M_i\oplus M_j$, and acts identically on $(M_i+M_j)/M_i$.}

\medskip
\proof
We will prove this only in the case when $i\neq j$ (as the case $i=j$ is even simpler). We know that $g_w$ acts identically on each $M_l$ with $l\in I\setminus\{j\}$ and on $(M_i\oplus M_j)/M_i$. This implies that each $M_l$ with $l\in I\setminus\{j\}$ is the maximal direct summand of $M$ such that all Newton polygon slopes of $(M_l,g_w\phi)$ are equal to $\alpha_l$ and that $M_i\oplus M_j$ is the maximal direct summand of $M$ such that all Newton polygon slopes of $(M_i\oplus M_j,g_w\phi)$ are equal to either $\alpha_i$ or $\alpha_j$. From this and the fact that $h_w\in \pmb{GL}_M(W(k))$ is an isomorphism between $(M,g_w\phi)$ and $(M,\phi)$, we get that $h_w$ leaves invariant each $M_l$ with $l\in I\setminus\{j\}$ as well as $M_i\oplus M_j$. Even more, from the second sentence of this proof we get that $h_w$ restricted to each $M_l$ with $l\in I\setminus\{j\}$ is an automorphism $h_{lw}$ of $(M_l,\phi)$ and moreover $h_w$ induces an automorphism of $((M_i\oplus M_j)/M_i,\phi)$ and thus an automorphism $h_{jw}$ of $(M_j,\phi)$. 

Let $h_{0w}:=\prod_{l\in I} h_{lw}\in\prod_{l\in I} \pmb{GL}_{M_l}(W(k))\leqslant \pmb{GL}_M(W(k))$; it is an automorphism of $(M,\phi)$. The element $\tilde h_w:=h_{0w}^{-1}h_w\in\pmb{GL}_M(W(k))$ has all the desired properties.\endproof

\medskip
To reach the desired contradiction we can assume that we have $h_w=\tilde h_w$, where $\tilde h_w$ is as in Fact 3.5.4. We first consider the case when $i\ne j$. 

\medskip\noindent
{\bf 3.5.5. The case $i\ne j$.} We assume that $i\neq j$ (i.e., $\alpha_j<\alpha_i$). As $h_w=\tilde h_w$, we have $u_w\in\Hom(M_j,M_i)$. From this and the relation $i\neq j$ we get that $u_w(w)=0$. As $\alpha_j<\alpha_i$, all Newton polygon slopes of $(\Hom(M_j,M_i),\phi)$ are positive. Therefore for each $V$ in $M_i$  the sequence $(\phi^m(p^{\ell-1}V\otimes y^*))_{m\ge 0}$ converges to $0$ and thus all the solutions of the equation $(7c)$ in $\Hom(M_j,M_i)[{1\over p}]\times M_i$ are of the form $(-\sum_{m=0}^{\infty} \phi^m(p^{\ell-1}V\otimes y^*),V)$. From this and the relation $u_w(w)=0$ we get the following identity
$$u_w=-\sum_{m=0}^{\infty} \phi^m(p^{\ell-1}w\otimes y^*).\leqno (7d)$$
We have the following two properties of the terms of the sum $(7d)$.
\medskip
{\bf (i)} All the terms of the sum of $(7d)$ belong to ${1\over p}\Hom(M_j,M_i)$ (this is so as $w$ and $y$ belong to a $W(k)$-basis for $M$ formed by elements of $\cup_{l\in I} M_l$ and therefore the element $\ell(w,y)$ can be defined as in Subsubsection 3.5.1 and it is equal to $\ell(x,y)=\ell$). Moreover, all but a finite number of these terms belong to $\Hom(M_j,M_i)$.  
\smallskip
{\bf (ii)} The term  $\phi^q(p^{\ell-1}w\otimes y^*)$ of the sum of $(7d)$ belongs to ${1\over p}\Hom(M_j,M_i)\setminus \Hom(M_j,M_i)$ (cf. property 3.5.2 (ii) and the fact that $a_{w,q}=a_{x,q}$).
\medskip
Let $\gamma$ be an invertible element of $W(k)$. Let $\bar\gamma\in k\setminus\{0\}$ be its reduction modulo $p$. Based on properties (i) and (ii), the condition that the element $u_{\gamma w}$ obtained as in $(7d)$ belongs to $\Hom(M_j,M_i)$ is expressed by $\bar\gamma$ being a solution of a system of polynomial equations in one variable which have coefficients in $k$ and which contain at least one polynomial of degree at least $p^q$. Therefore there exist such elements $\gamma$ with the property that we have $u_{\gamma w}\in {1\over p}\Hom(M_j,M_i)\setminus \Hom(M_j,M_i)$. Thus for such an element $\gamma$ we have $\gamma w\in \Lambda$ and $h_{\gamma w}=1_M+u_{\gamma w}\notin\pmb{GL}_M(W(k))$. Contradiction. 

\medskip\noindent
{\bf 3.5.6. Extra reduction steps.} To reach the desired contradiction we can assume that $i=j$ (i.e., $\alpha_i=\alpha_j$), cf. Subsection 3.5.5. As $h_w=\tilde h_w$ and $i=j$, to reach a contradiction we can assume based on Fact 3.5.4 (i) that $M_i=M$ (i.e., that $I=\{i\}$). Thus $(M,\phi)$ is isoclinic and we have $O=O_0=A_0\otimes_{\dbZ_p} W(k)$, cf. Subsection 1.2. 

\medskip\noindent
{\bf 3.5.7. Lemma.} {\it We recall that $v_w=w+u_w(w)$. Then we have $v_w\in\Lambda$.}

\medskip
\proof
Due to the definition of $q$, we have $\phi^s(g_w)\in\pmb{GL}_M(W(k))$ for all $s\in\{1,\ldots,q-1\}$ but $\phi^q(g_w)\notin\pmb{GL}_M(W(k))$. From this and the equation $(7a)$ we get that $\phi^s(h_w)\in\pmb{GL}_M(W(k))$ for all $s\in\{1,\ldots,q\}$ but $\phi^{q+1}(h_w)\notin\pmb{GL}_M(W(k))$. Thus we have $\phi^s(u_w)\in\End(M)$ for all $s\in\{1,\ldots,q\}$ but $\phi^{q+1}(u_w)\notin\End(M)$. 

Due to this and the identity $(7b)$ we get that $\phi^q(p^{\ell-1}v_w\otimes y^*)\notin\End(M)$. If  $\phi^q(p^{\ell-1}v_w\otimes y^*)\notin {1\over p}\End(M)$ or if $v_w\in pM$, then we have $\ell_0\ge \ell+1$ and this contradicts $(6a)$. Thus we have $\phi^q(p^{\ell-1}v_w\otimes y^*)\in {1\over p}\End(M)\setminus\End(M)$ and $v_w\in M\setminus pM$. Therefore $\phi^q(p^{\ell}v_w\otimes y^*)\in \End(M)\setminus p\End(M)$ and $v_w\in M\setminus pM$. But we also have $\phi^q(p^{\ell}x\otimes y^*)\in \End(M)\setminus p\End(M)$, cf. property 3.5.2 (ii). From  the last two sentences and the very definitions of $a_{w,q}$ and $a_{x,q}$, we get that $a_{v_w,q}=a_{x,q}$. From this and the relation $v_w\in M\setminus pM$ we conclude that $v_w\in\Lambda$.\endproof 

\medskip\noindent
{\bf 3.5.8. Lemma.} {\it Let $(u,v)\in \End(M)\times M$ be a solution of the equation $(7c)$. 

\medskip
{\bf (a)} Then we have $\{u, p^{\ell-1}v\otimes y^*\}\subset \End(M)\cap {1\over p}O$.

\smallskip
{\bf (b)} Let $v_1\in pM$. Then there exists a solution $(u+u_1,v+v_1)$ of the equation $(7c)$ with $u_1\in O$.}

\medskip
\proof
We have $pp^{\ell-1}v\otimes y^*\in p^{\ell}\End(M)\subseteq O$. It is easy to see that for each element $\tilde w\in O$, the equation $*+\tilde w=\phi(*)$ in $*$ has a solution in $\phi(O)=O=O_0$ and thus also in $\End(M)$. Let $\tilde u\in O$ be such that we have $\tilde u+pp^{\ell-1}v\otimes y^*=\phi(\tilde u)$. Thus $pu-\tilde u=\phi(pu-\tilde u)$ belongs to $A_0[{1\over p}]\cap\End(M)=A_0\subseteq O$. Therefore we have $u\in\End(M)\cap {1\over p}O$. Thus (a) holds. Part (b) follows from the fact that there exists $u_1\in O$ such that $u_1+p^{\ell}v_1\otimes y^* =\phi(u_1)$. \endproof

\medskip\noindent
{\bf 3.5.9. Morphisms between $k$-schemes.} Let $\scrM$ be the affine space (scheme) over $k$ defined naturally by the $k$-vector space $M/pM$. Let $\vph_q:\scrM\to\scrM$ be the morphism of $k$-schemes that takes $\bar *\in\scrM(k)=M/pM$ to the element of $\scrM(k)$ which is the reduction modulo $p$ of $p^{-a_{x,q}}\phi^q(*)\in M$, where $*\in M$ is an arbitrary lift of $\bar m$.  

The set $\im(\Lambda\to M/pM)$ is the set of $k$-valued points of the open, non-empty subscheme $\scrS:=\vph_q^{-1}(\scrM\setminus\{0\})$ of $\scrM$. For each solution $(u,v)\in \End(M)\times M$ of the equation $(7c)$, a similar argument to the proof of Lemma 3.5.8 shows that $v$ modulo $p$ determines $u$ modulo $A_0$ up to a finite number of possibilities. From this and the identity $w=v_w-u_w$, we get that the association that takes $(w,u_w)$ modulo $p$ to $v_w$ modulo $p$ has finite fibres. This association can be viewed as the one defined naturally (at the level of $k$-valued points) by a morphism of $k$-schemes whose codomain is $\scrS$ and whose domain has the same dimension $r$ as $\scrS$. By reasons of dimensions, we get that:
\medskip
{\bf (i)} There exists an open, non-empty subscheme $\scrV$ of $\scrS$ which has the property that each $k$-valued point $\bar v$ of $\scrV$ is of the form $v_w$ modulo $p$ for some elements $w\in\Lambda$ and $u_w\in\End(M)$ such that $(w,v_w):=(w,w+u_w(w))$ is a solution of the equation $(7c)$. 
\medskip
Let $\bar O:={1\over p}O/O$ and let $\bar E$ be the image of $\End(M)\cap {1\over p}O$ in $\bar O$. Both $\bar O$ and $\bar E$ are $k$-vector spaces. Let $\scrO$ and $\scrE$ be the affine spaces (schemes) over $k$ defined naturally by the $k$-vector spaces $\bar O$ and $\bar E$ (respectively). Let $\bar\phi:\scrO\to\scrO$ be the morphism which takes a $k$-valued point of $\scrO$ defined by some element $o\in {1\over p}O$ to the $k$-valued point of $\scrO$ defined by the element $\phi(o)-o\in {1\over p}O$ (we think of $\bar\phi$ as a finite, surjective endomorphism of $\dbG_a^{r^2}$). Let $\scrF:=\scrE\cap \bar\phi^{-1}(\scrE)$. Thus $\scrF$ is a closed subscheme of $\scrO$ equipped with a morphism $m_1:\scrF\to\scrE$ induced from $\bar\phi$ (we think of $m_1$ as a homomorphism between closed subgroup schemes of $\dbG_a^{r^2}$). Based on Lemma 3.5.8 (a) we can speak about the natural images $\bar u_w$ and $\bar v_{w,y}$ of $u_w$ and $p^{\ell-1}v_w\otimes y^*$ (respectively) in $\bar O$ and thus about $k$-valued points (denoted in the same way) $\bar u_w\in\scrF(k)$ and $\bar v_{w,y}\in\scrE(k)$ with the property that $m_1$ maps $\bar u_w$ to $\bar v_{w,y}$. 

We have a natural morphism of $k$-schemes $m_2:\scrV\to \scrE$ which at the level of $k$-valued points maps a $k$-valued point of $\scrV$ represented by an element $v\in\Lambda$ to the $k$-valued point of $\scrE$ defined by the image of $p^{\ell-1}v\otimes y^*\in p^{\ell-1}\End(M)\subseteq {1\over p}O$ in $\bar O$. From the property (i) and the previous paragraph we get that the natural morphism
$$\iota:\scrV\times_{\scrE} \scrF\to\scrV$$ 
associated to the fibre product of $m_1$ and $m_2$, is surjective. As the morphism $\bar\phi:\scrO\to\scrO$ has finite fibres and it is of finite type, the finite type morphism $\iota$ is quasi-finite and therefore it is generically finite. From this, the property (i), and Lemma 3.5.8 (b) we get that:

\medskip
{\bf (ii)} There exists an element $v\in\Lambda$ that defines naturally a $k$-valued point of $\scrV$ and there exists a finite subset $\bar\Gamma$ of $k$ such that for each algebraically closed field $k_1$ that contains $k$ and for every invertible element $\gamma$ of $W(k_1)$ whose reduction modulo $p$ does not belong to $\bar\Gamma$, the following equation in $U$ 
$$U+p^{\ell-1}\gamma v\otimes y^*=(\phi\otimes \sigma_{k_1})(U)\leqno (7e)$$
obtained from $(7c)$ by replacing $(V,\phi)$ with $(\gamma v,\phi\otimes\sigma_{k_1})$, possesses a solution in $\End(M\otimes_{W(k)} W(k_1))$. Here $\sigma_{k_1}$ is the Frobenius automorphism of the ring $W(k_1)$ of Witt vectors with coefficients in $k_1$. 

\medskip\noindent
{\bf 3.5.10. Good choice of $\gamma$.} We will take $k_1$ to be an algebraic closure of $k((X))$, where $X$ is an independent variable. We identify $W(k)[[X]]$ with a $W(k)$-subalgebra of $W(k_1)$ that contains the invertible element $X=(X,0,\ldots)$ of $W(k_1)$. We will take $\gamma:=\tau X$, where $\tau$ is an invertible element of $W(k)$. We have $\sigma_{k_1}(X)=X^p$. For this choice of $\gamma$, the equation $(7e)$ has (up to addition of elements in the free $\dbZ_p$-module ${1\over p}A_0$ of rank $r^2$) a unique solution
$$u_{\tau X}=-\sum_{m=0}^{\infty} X^{p^m}\phi^m(p^{\ell-1}\tau v\otimes y^*)\leqno (7f)$$
in ${1\over p}O\otimes_{W(k)} W(k_1)$. In fact we have $u_{\tau X}\in {1\over p}\End(M)\otimes_{W(k)} W(k)[[X]]$. As $\tau v\in\Lambda$, from the property 3.5.2 (ii) we get that the term  $X^{p^q}\phi^q(p^{\ell-1}\tau v\otimes y^*)$ of $(7f)$ does not belong to $\End(M)\otimes_{W(k)} W(k_1)$. The last two sentences imply that the intersection $(u_{\tau X}+{1\over p}A_0)\cap [\End(M)\otimes_{W(k)} W(k_1)]$ is empty and therefore we reached the desired contradiction.  

\medskip\noindent
{\bf 3.5.11. End of the proof.} The contradiction we reached implies that $n_{\pmb{GL}_M}>\ell-1$. Thus $n_{\pmb{GL}_M}=\ell=\ell_{\pmb{GL}_M}$. This ends the proof of Theorem 1.3 (b) and therefore also of the Main Theorem A.\endproof

\medskip\smallskip\noindent
{\bf 3.6. Remarks.} Suppose $(M,\phi)$ is a direct sum of isoclinic latticed $F$-isocrystals over $k$. 

\medskip
{\bf (a)} We have a direct sum decomposition $\grg=(\grg\cap L_+)\oplus(\grg\cap L_0)\oplus(\grg\cap L_-)$ of $W(k)$-modules. Thus $O_G=\grg\cap O$ and therefore $\grg/O_G\subseteq \End(M)/O$. From this inclusion we easily get the following {\it monotony properties}: we have $\ell_G\le\ell_{\pmb{GL}_M}$ and therefore (cf. Main Theorem A) we also have $n_G\le n_{\pmb{GL}_M}$.

\smallskip
{\bf (b)} We assume that $(M,\phi)$ is the Dieudonn\'e module of $D$; thus $\ell_{\pmb{GL}_M}=\ell_D$. We will use the notations of Subsection 3.5. We also assume that there exist elements $x,y\in\scrB$ such that $\ell(x,y)\ge 2$ and we have $x\in\scrB_i$ and $y\in\scrB_j$ with $\alpha_j<\alpha_i$. Let $g\in\pmb{GL}_M(W(k))$ be the element that fixes each $z\in\scrB\setminus\{y\}$ and that takes $y$ to $y+p^{\ell(x,y)-1}x$. Let $D_g$ be a $p$-divisible group over $k$ whose Dieudonn\'e module is isomorphic to $(M,g\phi)$. Then $D_g[p^{\ell(x,y)-1}]$ is isomorphic to $D[p^{\ell(x,y)-1}]$ and $D_g$ has the same Newton polygon as $D$ (as $\alpha_j\neq\alpha_i$). If  by chance we also have an identity $\ell_D=\ell(x,y)$, then  Subsection 3.5 can be easily adapted to give us that, up to a replacement of $x\in\scrB_i$ by a multiple of it with an invertible element of $W(k)$, we can assume that $D_g$ is not isomorphic to $D$.

\medskip\smallskip\noindent
{\bf 3.7. Variant of 1.3 (b).} Let $(M,\phi,G)$ be a latticed $F$-isocrystal with a group over $k$ such that Assumption 1.1.1 holds. We assume that the following two conditions hold:

\medskip
{\bf (i)} we have $n_G\ge 1$ and a direct sum decomposition $\grg=(\grg\cap L_+)\oplus (\grg\cap L_0)\oplus (\grg\cap L_-)$ (or $\grg=(\grg\cap L_+)\oplus [\grg\cap (L_0\oplus L_-)]$ or $\grg=[\grg\cap (L_+\oplus L_0)]\oplus (\grg\cap L_-)$) of $W(k)$-modules;

\smallskip
{\bf (ii)} for all $q\in\dbN$, there exists a $W(k)$-basis $\scrB$ for $M$ and a sequence of integers $(a_{z,q})_{z\in\scrB}$ such that certain subsets of $\{x\otimes y^*|x,y\in\scrB\}$ are $W(k)$-bases for all direct summands of $\grg$ listed in (i) and moreover we have $M=\oplus_{z\in\scrB} W(k)p^{-a_{z,q}}\phi^q(z)$.

\medskip
Then the proof of Theorem 1.3 (b) (see Subsection 3.5) can be entirely adapted to give us that $n_G=\ell_G$. We only add here two things. First, if by chance in Subsubsection 3.5.2 we have $\ell(x,y)=\ell_{\pmb{GL}_M}$ with $x\in\scrB_i$ and $y\in\scrB_j$ such that $\alpha_j>\alpha_i$, then one needs to use $\phi^{-q}$ (instead of $\phi^{q}$) with $q\in\dbN$ in order to reach the desired contradiction. Second, if we have $\grg=(\grg\cap L_+)\oplus [\grg\cap (L_0\oplus L_-)]$ (resp. $\grg=[\grg\cap (L_+\oplus L_0)]\oplus (\grg\cap L_-)$), then one needs to use Lemma 2.4 in order to be able to treat $\grg\cap (L_0\oplus L_-)$ (resp. $\grg\cap (L_+\oplus L_0)$) in the same manner as $\grg\cap L_-$ (resp. as $\grg\cap L_+$).

\bigskip
\noindent
{\boldsectionfont 4. Direct applications to $p$-divisible groups} 
\bigskip
 
In this Section we prove the results stated in Subsubsections 1.4.2 to 1.4.4 (see Subsections 4.5 to 4.7). In Subsections 4.1 to 4.4 we introduce basis invariants of $p$-divisible groups over $k$ and we present basic properties of them that are needed in Subsections 4.5 to 4.7. Until the end we will assume that $(M,\phi)$ is the Dieudonn\'e module of $D$.

\medskip\smallskip\noindent
{\bf 4.1. Definitions.} {\bf (a)} Let $q\in\dbN$. Let $\alpha_D(q)\in\dbN\cup\{0\}$ be the greatest number such that we have $\phi^q(M)\subseteq p^{\alpha_D(q)}M$. Let $\beta_D(q)\in\dbN\cup\{0\}$ be the smallest number such that we have $p^{\beta_D(q)}M\subseteq \phi^q(M)$. Let $\delta_D(q):=\beta_D(q)-\alpha_D(q)$; as $p^{\beta_D(q)}M\subseteq \phi^q(M)\subseteq p^{\alpha_D(q)}M$, we have $\delta_D(q)\in\dbN\cup\{0\}$. 
 
\smallskip
{\bf (b)} We assume that $D$ is isoclinic. Let $m:=g.c.d.\{c,d\}$. Let $(c_1,d_1,r_1):=({c\over m},{d\over m},{r\over m})$. Let $u_D:=\text{sup}\{0,\beta_D(r_1n)-d_1n|n\in\dbN\}$. Let $\tilde u_D:=\text{sup}\{0,d_1n-\alpha_D(r_1n)|n\in\dbN\}$. Let $v_D:=\text{sup}\{0,\beta_D(rn)-dn|n\in\dbN\}$. Let $\tilde v_D:=\text{sup}\{0,dn-\alpha_D(rn)|n\in\dbN\}$. Proposition 4.3 (c) and (b) below will imply that $u_D$, $\tilde u_d$, $v_D$, and $\tilde v_D$ are non-negative integers and therefore that in their definition we can replace $\text{sup}$ by $\max$.

\smallskip
{\bf (c)} Let $(M_1,\phi_1)$ and $(M_2,\phi_2)$ be the Dieudonn\'e modules of two isoclinic $p$-divisible groups $D_1$ and $D_2$ (respectively) over $k$. Let $\alpha_1$ and $\alpha_2$ be the unique Newton polygon slopes of $D_1$ and $D_2$ (respectively). Let $\ell_{D_1,D_2}\in\dbN\cup\{0\}$ be the smallest number that has the following property: 

\medskip
{\bf (i)} if $\alpha_1\le\alpha_2$, then for all $q\in\dbN$ the $W(k)$-module 
$\phi^q(p^{\ell_{D_1,D_2}}\Hom(M_1,M_2))=p^{\ell_{D_1,D_2}}\Hom(\phi^q(M_1),\phi^q(M_2))$ is included in $\Hom(M_1,M_2)$;

\smallskip
{\bf (ii)} if $\alpha_1>\alpha_2$, then  for all $q\in\dbN$ the $W(k)$-module 
$\phi^{-q}(p^{\ell_{D_1,D_2}}\Hom(M_1,M_2))=p^{\ell_{D_1,D_2}}\Hom(\phi^{-q}(M_1),\phi^{-q}(M_2))$ is included in $\Hom(M_1,M_2)$.

\medskip
{\bf (d)} Let $\eps_D:=\eps_{\pmb{GL}_M}$, where the number $\eps_{\pmb{GL}_M}\in\{0,1\}$ is as in Scholium 3.5.1. 

\medskip
\medskip\noindent
{\bf 4.1.1. Remark.} If $cd=0$, then $O_0=\End(M)$ and therefore $\ell_D=\eps_D=0$. If $c,d\ge 1$ and $D$ is ordinary (i.e., isomorphic to $(\dbQ_p/\dbZ_p)^c\oplus (\pmb{\mu}_{p^{\infty}})^d$), then $O=\End(M)$ and the two-sided ideal of the $W(k)$-algebra $\End(M)$ generated by $O_+\oplus O_-$ is $\End(M)$; thus $\ell_D=\eps_D=1$. If $c,d\ge 1$ and $D$ is not isomorphic to $(\dbQ_p/\dbZ_p)^c\oplus (\pmb{\mu}_{p^{\infty}})^d$, then $\End(M)\neq O$ and therefore $\ell_D>0$ and $\eps_D=0$; moreover $\ell_D\in\dbN$ is the smallest number such that we have $p^{\ell_D}\End(M)\subseteq O$ (cf. rule 1.2 (b)).

\medskip\smallskip\noindent
{\bf 4.2. Simple properties.} In this Subsection we list few simple properties of the invariants we have introduced so far.

\medskip\noindent
{\bf 4.2.1. Fact.} {\it We have $n_D=n_{D^{\text{t}}}$ and $\ell_D=\ell_{D^{\text{t}}}$.}

\medskip
\proof
We show that $n_D\le n_{D^{\text{t}}}$. Let $C$ be a $p$-divisible group of codimension $c$ and dimension $d$ over $k$. If $C[p^{n_{D^{\text{t}}}}]$ is isomorphic to $D[p^{n_{D^{\text{t}}}}]$, then taking Cartier duals we get that $C^{\text{t}}[p^{n_{D^{\text{t}}}}]$ is isomorphic to $D^{\text{t}}[p^{n_{D^{\text{t}}}}]$ and thus that $C^{\text{t}}$ is isomorphic to $D^{\text{t}}$. Taking Cartier duals, we get that $C$ is isomorphic to $D$. This implies that $n_D\le n_{D^{\text{t}}}$. As $D$ is the Cartier dual of  $D^{\text{t}}$, we also have $n_{D^{\text{t}}}\le n_D$. Thus $n_D=n_{D^{\text{t}}}$.

As $(M^*,p\phi)$ is the Dieudonn\'e module of $D^{\text{t}}$, under the natural identification $\End(M^*)=\End(M)$, the level module of $(M^*,p\phi)$ gets identified with $O$. Thus we have $\ell_D=\ell_{D^{\text{t}}}$.\endproof

\medskip\noindent
{\bf 4.2.2. Fact.} {\it The following three properties hold:

\medskip
{\bf (a)} for all $q\in\dbN$, we have inclusions 
$$M\subseteq p^{-\beta_D(q)}\phi^q(M)\subseteq p^{-\delta_D(q)}M\leqno (8a)$$
which are optimal in the sense that we also have
$$M\nsubseteq p^{-\beta_D(q)+1}\phi^q(M)\;\;\text{and}\;\;p^{-\beta_D(q)}\phi^q(M)\nsubseteq p^{-\delta_D(q)+1}M;\leqno (8b)$$

\indent
{\bf (b)} for all $q\in\dbN$, we have $\alpha_{D^{\text{t}}}(q)=q-\beta_D(q)$ and $\beta_{D^{\text{t}}}(q)=q-\alpha_D(q)$;

\smallskip
{\bf (c)} if $D$ is isoclinic, then we have $u_D=\tilde u_{D^{\text{t}}}$, $\tilde u_D=u_{D^{\text{t}}}$, $v_D=\tilde v_{D^{\text{t}}}$, and $\tilde v_D=v_{D^{\text{t}}}$.}

\medskip
\proof
Part (a) follows from the very Definition 4.1 (a). As $p^{\beta_D(q)}M\subseteq \phi^q(M)\subseteq p^{\alpha_D(q)}M$, we have $p^{-\alpha_D(q)}M^*\subseteq \phi^q(M^*)\subseteq p^{-\beta_D(q)}M^*$ i.e., $p^{q-\alpha_D(q)}M^*\subseteq (p\phi)^q(M^*)\subseteq p^{q-\beta_D(q)}M^*$. As $(M^*,p\phi)$ is the Dieudonn\'e module of $D^{\text{t}}$ and due to $(8a)$ and $(8b)$, we get that (b) holds. We prove (c). Due to (b) we have an equality $\beta_D(r_1n)-d_1n=r_1n-\alpha_{D^{\text{t}}}(r_1n)-d_1n=c_1n-\alpha_{D^{\text{t}}}(r_1n)$. This implies that $u_D=\tilde u_{D^{\text{t}}}$. By replacing $D$ with $D^{\text{t}}$, we get that $u_{D^{\text{t}}}=\tilde u_D$. Similar arguments show that $v_D=\tilde v_{D^{\text{t}}}$ and $\tilde v_D=v_{D^{\text{t}}}$. Thus (c) holds.\endproof

\medskip\noindent
{\bf 4.2.3. Lemma.} {\it We assume that $D$ is isoclinic. Let $\alpha:={d\over r}\in\dbQ\cap [0,1]$ be its unique Newton polygon slope. Then the following two properties hold:

\medskip
{\bf (a)} we have $\alpha_D(q)\le q\alpha\le \beta_D(q)$;

\smallskip
{\bf (b)} if $\alpha_D(q)=q\alpha$ (or if $\beta_D(q)=q\alpha$), then we have $\alpha_D(q)=\beta_D(q)=q\alpha$.}

\medskip
\proof
As $\vph_q:=p^{-\alpha_D(q)}\phi^q$ is a $\sigma^q$-linear endomorphism of $M$, the Newton polygon slopes of $\vph_q$ are on one hand non-negative and on the other hand are all equal to $q\alpha-\alpha_D(q)$. Thus $\alpha_D(q)\le q\alpha$. If $q\alpha=\alpha_D(q)$, then all the Newton polygon slopes of $\vph_q:M\to M$ are $0$ and therefore we have $\vph_q(M)=M$. This implies that $\beta_D(q)=\alpha_D(q)=q\alpha$. The part involving $\beta_D(q)$ is proved in the same way but working with $p^{\beta_D(q)}\phi^{-q}$.\endproof 

\medskip\noindent
{\bf 4.2.4. Lemma.} {\it {\bf (a)} If $D$ is either $F$-cyclic or special, then $D$ is also quasi-special.

\smallskip
{\bf (b)} The class $\scrQ_{c,d}$ of isomorphism classes of quasi-special $p$-divisible groups of codimension $c$ and dimension $d$ over $k$, is a finite set.}

\medskip
\proof
Each isoclinic special $p$-divisible group over $k$ is isoclinic quasi-special. Each $F$-cyclic $p$-divisible group over $k$ is a direct sum of $F$-circular $p$-divisible groups over $k$. Based on the last two sentences, it suffices to prove (a) in the case when $D$ is $F$-circular. Let $\pi$ be an $r$-cycle of $J_r$ such that $D$ is isomorphic to $C_{\pi}$. We have $\phi_{\pi}^r(M)=p^dM$ and therefore $C_{\pi}$ is isoclinic quasi-special of Newton polygon slope ${d\over r}$. Thus (a) holds.

To prove (b) it suffices to show that for all pairs $(c,d)\in (\dbN\cup\{0\})^2$ with $c+d>0$, the class $\scrI_{c,d}$ of isomorphism classes of isoclinic quasi-special $p$-divisible groups of codimension $c$ and dimension $d$ over $k$, is a finite set. We assume that $D$ is isoclinic quasi-special. Then we have $\phi^r(M)=p^dM$. Therefore $\vph_r:=p^{-d}\phi^r:M\to M$ is a $\sigma^r$-linear automorphism of $M$. Let $M_{W(\dbF_{p^r})}:=\{x\in M|\vph_r(x)=x\}$. We have $M_{W(\dbF_{p^r})}\otimes_{W(\dbF_{p^r})} W(k)=M$. Moreover $\phi(M_{W(\dbF_{p^r})})\subseteq M_{W(\dbF_{p^r})}$. Therefore the Dieudonn\'e module $(M,\phi)$ is definable over the finite field $\dbF_{p^r}$. Thus every isoclinic quasi-special $p$-divisible group of codimension $c$ and dimension $d$ over $k$ has a Dieudonn\'e module over $k$ which: (i) is isomorphic to $(M,g\phi)$ for a suitable element $g\in\pmb{GL}_M(W(k))$, and (ii) it is definable over $\dbF_{p^r}$.

Let $M=F^1\oplus F^0$ be a direct decomposition such that $F^1/pF^1$ is the kernel of $\phi$ modulo $p$. We have $\phi({1\over p}F^1+F^0)=M$. Thus the cocharacter $\mu:\dbG_m\to\pmb{GL}_M$ that fixes $F^0$ and that acts on $F^1$ via the inverse of the identity character of $\dbG_m$, is a Hodge cocharacter of $(M,\phi,\pmb{GL}_M)$ in the sense of [Va1, Subsubsect. 2.2.1 (d)]. Thus the triple $(M,\phi,\pmb{GL}_M)$ is a latticed $F$-isocrystal with a group over $k$ for which the $W$-condition of loc. cit. holds. From the Atlas Principle applied to $(M,\phi,\pmb{GL}_M)$ and to an emphasized family of tensors indexed by the empty set (see [Va1, Thm. 5.2.3]), we get that the set of isomorphism classes of Dieudonn\'e modules over $k$ which are of the form $(M,g\phi)$ with $g\in\pmb{GL}_M(W(k))$ and which are definable over the finite field $\dbF_{p^r}$ is finite. From this and the classical Dieudonn\'e theory, we get that the class $\scrI_{c,d}$ is a finite set.\endproof

\medskip\smallskip\noindent
{\bf 4.3. Proposition.}  {\it We assume that $D$ is isoclinic. Then the following six properties hold:

\medskip
{\bf (a)} we have $\ell_D=\max\{\delta_D(q)|q\in\dbN\}$;

\smallskip
{\bf (b)} if $\alpha:={d\over r}$, then we have $\lim_{q\to\infty} {{\beta_D(q)\over q}}=\lim_{q\to\infty} {{\alpha_D(q)\over q}}=\alpha$;

\smallskip
{\bf (c)} if $M_0$ (resp. $\tilde M_0$) is the $W(k)$-submodule of $M$ generated by elements fixed by $\vph_{r_1}:=p^{-d_1}\phi^{r_1}$ (resp. by $\vph_r:=p^{-d}\phi^{r}$), then $u_D$ (resp. $v_D$) is finite and it is the smallest non-negative integer such that $p^{u_D}$ (resp. $p^{v_D}$) annihilates $M/M_0$ (resp. $M/\tilde M_0$);

\smallskip
{\bf (d)} if $M_1$ (resp. $\tilde M_1$) is the smallest $W(k)$-submodule of $M[{1\over p}]$ which is generated by elements fixed by $\vph_{r_1}$ (resp. by $\vph_r$) and which contains $M$, then $\tilde u_D$ (resp. $\tilde v_D$) is finite and it is the smallest non-negative integer such that $p^{\tilde u_D}$ (resp. $p^{\tilde v_D}$) annihilates $M_1/M$ (resp. $\tilde M_1/M$);

\smallskip
{\bf (e)} we have $u_D=\tilde u_D$ (resp. $v_D=\tilde v_D$);

\smallskip
{\bf (f)} we have $u_D\le\ell_D$.}

\medskip
\proof
We prove (a). The $W(k)$-span of endomorphisms of $(M,\phi)$ is $O=O_0$. The number $\ell_D$ is the smallest number such that $p^{\ell_D}\End(M)\subseteq O\subseteq \End(M)$, cf. Example 1.2.2. As
$O=\cap_{q\in\dbN\cup\{0\}} \phi^q(\End(M))=\cap_{q\in\dbN\cup\{0\}} \End(\phi^q(M))=\cap_{q\in\dbN\cup\{0\}} \End(p^{-\beta_D(q)}\phi^q(M))$, $\ell_D$ is the smallest (non-negative) integer such that we have $p^{\ell_D}\End(M)\subseteq \End(p^{-\beta_D(q)}\phi^q(M))$ for all $q\in\dbN$. Thus from $(8a)$ and $(8b)$ we get that $\ell_D$ is the smallest integer which is greater or equal to $\delta_D(q)$ for all $q\in\dbN$. From this (a) follows.

We prove (b). From (a) we get that $\delta_D(q)=\beta_D(q)-\alpha_D(q)\le\ell_D$. Thus ${{\alpha_D}\over q}\le \alpha\le {{\beta_D}\over q}\le {{\alpha_D}\over q}+{{\ell_D}\over q}$, cf. Lemma 4.2.3 (a). From these inequalities we get that (b) holds. 

We will prove (c) only for $M_0$ as the case of $\tilde M_0$ is argued in the same manner. We have $M_0=\cap_{n\in\dbN\cup\{0\}} \vph_{r_1}^n(M)$. As $\vph_{r_1}^n(M)=p^{-d_1n}\phi^{r_1n}(M)$, from $(8a)$ and $(8b)$ we get that $p^{\beta_D(r_1n)-d_1n}M\subseteq\vph_{r_1}^n(M)$ and $p^{\beta_D(r_1n)-d_1n-1}M\nsubseteq\vph_{r_1}^n(M)$. Thus the smallest non-negative number $s$ such that we have $p^sM\subseteq\vph_{r_1}^n(M)$ for all $n\in\dbN$ is $\text{sup}\{0,\beta_D(r_1n)-d_1n|n\in\dbN\}$ and therefore it is $u_D$. Thus (c) holds.

We will prove (d) only for $M_1$ as the case of $\tilde M_1$ is argued in the same manner. The $W(k)$-submodule $M_1^*$ of $M^*$ is the largest $W(k)$-submodule of $M^*$ generated by elements fixed by $p^{-d_1}\phi^{r_1}$. As $(M^*,p\phi)$ is the Dieudonn\'e module of $D^{\text{t}}$, the analogue of $p^{-d_1}\phi^{r_1}$ for $D^{\text{t}}$ is the $\sigma^{r_1}$-linear automorphism $p^{-c_1}(p\phi)^{r_1}=p^{-d_1}\phi^{r_1}$ of $M^*[{1\over p}]$. Thus from (c) applied to $D^{\text{t}}$, we get that $u_{D^{\text{t}}}$ is the smallest non-negative integer  with the property that $p^{u_{D^{\text{t}}}}$ annihilates $M^*/M_1^*$. As $\tilde u_D=u_{D^{\text{t}}}$ (see Fact 4.2.2 (c)) and as the $W(k)$-modules $M_1/M$ and $M^*/M_1^*$ are isomorphic, we get that $\tilde u_D$ is the smallest non-negative integer such that $p^{\tilde u_D}$ annihilates $M_1/M$. Thus (d) holds.

We will prove (e) for $u_D$ and $\tilde u_D$ as the case of $v_D$ and $\tilde v_D$ is argued in the same manner. As $p^{\tilde u_D}M_1\subseteq M$ and as $p^{\tilde u_D}M_1$ is $W(k)$-generated by elements fixed by $\vph_{r_1}$, we have $p^{\tilde u_D}M_1\subseteq M_0$. Thus $p^{\tilde u_D}$ annihilates $M/M_0$ and therefore $\tilde u_D\ge u_D$, cf. (c). A similar argument shows that $M_1\subseteq p^{-u_D}M_0$ and that $\tilde u_D\le u_D$. Thus $u_D=\tilde u_D$ i.e., (e) holds. 

We prove (f). Each endomorphism of $(M,\phi)$ maps $M_0$ to $M_0$. Thus $O=O_0\subseteq\End(M_0)$. But due to (c), the smallest number $s\in\dbN\cup\{0\}$ such that $p^s\End(M)\subseteq\End(M_0)$ is $u_D$. As $p^{\ell_D}\End(M)\subseteq O\subseteq\End(M_0)$, we get that $u_D\le\ell_D$. Thus (f) holds.\endproof

\medskip\noindent
{\bf 4.3.1. Remark.} We have $M_0\subseteq \tilde M_0\subseteq M$ and $M\subseteq \tilde M_1\subseteq M_1$. Thus $v_D=\tilde v_D\le u_D=\tilde u_D$. The $W(k)$-submodule $M_0$ (resp. $\tilde M_0$) of $M$ is the largest one with the property that $(M_0,\phi)$ (resp. $(\tilde M_0,\phi)$) is the Dieudonn\'e module of an isoclinic special (resp. isoclinic quasi-special) $p$-divisible group over $k$. Thus we call $u_D=\tilde u_D$ (resp. $v_D=\tilde v_D$) the {\it Manin height} (resp. the {\it Manin quasi-height}) of $D$, cf. [Ma, Ch. III, Sect. 2]. Similarly, the $W(k)$-submodule $M_1$ (resp. $\tilde M_1$) of $M[{1\over p}]$ is the smallest one with the properties that it contains $M$ and that $(M_1,\phi)$ (resp. $(\tilde M_1,\phi)$) is the Dieudonn\'e module of an isoclinic special (resp. isoclinic quasi-special) $p$-divisible group over $k$.

\medskip\smallskip\noindent
{\bf 4.4. Proposition.} {\it We assume that $D=\prod_{i\in I} D_i$ is a product of at least two non-trivial isoclinic $p$-divisible groups over $k$. Then the following three properties hold:

\medskip
{\bf (a)} for all $i,j\in I$ with $i\neq j$, we have $\ell_{D_i,D_j}=\ell_{D_j,D_i}$;

\smallskip
{\bf (b)} we have $\ell_{D_i}=\ell_{D_i,D_i}$ and $\ell_D=\max\{\eps_D,\ell_{D_i},\ell_{D_i,D_j}|i\in I,j\in I\setminus\{i\}\}$;

\smallskip
{\bf (c)} if $i,j\in I$ with $i\neq j$ and if the Newton polygon slope $\alpha_i$ of $D_i$ is less or equal to the Newton polygon slope $\alpha_j$ of $D_j$, then we have $\ell_{D_i,D_j}=\max\{0,\beta_{D_i}(q)-\alpha_{D_j}(q)|q\in\dbN\}$.}

\medskip
\proof
Let $M=\bigoplus_{i\in I} M_i$ be the direct sum decomposition such that $(M_i,\phi)$ is the Dieudonn\'e module of $D_i$. Let $\scrB_i$ be a $W(k)$-basis for $M_i$. Let $\scrB:=\cup_{i\in I} \scrB_i$. For $x,y\in\scrB$, let $\ell(x,y)\in\dbN\cup\{0\}$ be defined as in Scholium 3.5.1. We have $\ell_{D_i,D_j}=\max\{\ell(y,x)|x\in\scrB_i,y\in\scrB_j\}$, cf. the very definitions. Thus (a) is a particular case of Formula $(6b)$. As $D_i$ is isoclinic, we have $\eps_{D_i}=0$. Thus (b) is a particular case of Formulas $(6a)$, $(5a)$, $(5b)$, and $(5c)$. We prove (c). Due  to $(8a)$, for all $q\in\dbN$ we have
$$\Hom(\phi^q(M_i),\phi^q(M_j))\subseteq \Hom(p^{\beta_{D_i}(q)}M_i,p^{\alpha_{D_j}(q)}M_j)=p^{\alpha_{D_j}(q)-\beta_{D_i}(q)}\Hom(M_i,M_j).$$
From this and making use of $(8a)$ and $(8b)$ we get that there exist a direct summand of $\Hom(\phi^q(M_i),\phi^q(M_j))$ which is also a direct summand of $p^{\alpha_{D_j}(q)-\beta_{D_i}(q)}\Hom(M_i,M_j)$. Thus the smallest number $s\in\dbN\cup\{0\}$ with the property that for all $q\in\dbN$ the $W(k)$-module $p^s\phi^q(\Hom(M_i,M_j))$ is included in $\Hom(M_i,M_j)$, is $\max\{0,\beta_{D_i}(q)-\alpha_{D_j}(q)|q\in\dbN\}$. From this and the rule (i) of the Definition 4.1 (c), we get that (c) holds.\endproof

\medskip\noindent
{\bf 4.4.1. Example.} We assume that $D$ is isoclinic and that $d< r\le 2d$. Thus $\alpha:={d\over r}\in\dbQ\cap [{1\over 2},1)$. For $q\in\dbN$ we have $\beta_{D^{\text{t}}}(q)-\alpha_D(q)=q-2\alpha_D(q)$, cf. Fact 4.2.2 (b). From this and Proposition 4.4 (c) we get that $\ell_{D^{\text{t}},D}=\max\{0,q-2\alpha_D(q)|q\in\dbN\}$. As $\eps_D=0$, from Proposition 4.4 (a) and (b) we get that $\ell_{D\oplus D^{\text{t}}}=\max\{\ell_D,\ell_{D^{\text{t}}},\ell_{D^{\text{t}},D}\}$. As $\ell_D=\ell_{D^{\text{t}}}=\max\{\delta_D(q)|q\in\dbN\}$ (cf. Fact 4.2.1 and Proposition 4.3 (a)), we conclude that 
$$\ell_{D\oplus D^{\text{t}}}=\max\{\delta_D(q),q-2\alpha_D(q)|q\in\dbN\}.$$

\noindent
{\bf 4.5. Proof of 1.4.3.} We assume that $D=\prod_{i\in I} D_i$ is a product of at least two non-trivial isoclinic $p$-divisible groups over $k$. Let $M=\bigoplus_{i\in I} M_i$ be the direct sum decomposition defined by the product decomposition $D=\prod_{i\in I} D_i$. As $n_D=\ell_D$ and as for $i\in I$ we have $n_{D_i}=\ell_{D_i}$ (cf. Corollary 1.4.2), based on Proposition 4.4 (a) and (b), to prove Proposition 1.4.3 it suffices to show that for all $i,j\in I$ with $i\neq j$ we have
$$\ell_{D_i,D_j}\le\max\{0,\ell_{D_i}+\ell_{D_j}-1\}.\leqno (9)$$
As $\ell_{D_i,D_j}=\ell_{D_j,D_i}$ (see Proposition 4.4 (a)), to check the inequality (9) we can assume that $\alpha_i\le\alpha_j$. We have $\alpha_{D_i}(q)\le q\alpha_i\le q\alpha_j\le\beta_{D_j}(q)$, cf. Lemma 4.2.3 (a). Thus $\alpha_{D_i}(q)\le\beta_{D_j}(q)$. Based on Proposition 4.4 (c), to prove the inequality (9) it suffices to show that for all $q\in\dbN$ we have $\beta_{D_i}(q)-\alpha_{D_j}(q)\le\max\{0,\ell_{D_i}+\ell_{D_j}-1\}$. We have $\delta_{D_i}(q)+\delta_{D_j}(q)\le\ell_{D_i}+\ell_{D_j}$, cf. Proposition 4.3 (a). From this and the inequality $\alpha_{D_i}(q)\le \beta_{D_j}(q)$ we get:
$$\beta_{D_i}(q)-\alpha_{D_j}(q)=\delta_{D_i}(q)+\alpha_{D_i}(q)+\delta_{D_j}(q)-\beta_{D_j}(q)\le\ell_{D_i}+\ell_{D_j}+\alpha_{D_i}(q)-\beta_{D_j}(q)\le\ell_{D_i}+\ell_{D_j}.\leqno (10)$$
If we have an equality $\beta_{D_i}(q)-\alpha_{D_j}(q)=\ell_{D_i}+\ell_{D_j}$, then $\beta_{D_j}(q)=\alpha_{D_i}(q)=q\alpha_i=q\alpha_j$ and therefore also $\beta_{D_i}(q)=q\alpha_i$ and $\alpha_{D_j}(q)=q\alpha_j=q\alpha_i$ (cf. Lemma 4.2.3 (b)). Thus the assumption that $\beta_{D_i}(q)-\alpha_{D_j}(q)=\ell_{D_i}+\ell_{D_j}$ implies that $\beta_{D_i}(q)-\alpha_{D_j}(q)=0=\ell_{D_i}+\ell_{D_j}\le\max\{0,\ell_{D_i}+\ell_{D_j}-1\}$. From this and (10) we get that the inequality $\beta_{D_i}(q)-\alpha_{D_j}(q)\le\max\{0,\ell_{D_i}+\ell_{D_j}-1\}$ always holds; therefore the inequality (9) holds. This ends the proof of Proposition 1.4.3.\endproof 

\medskip\smallskip\noindent
{\bf 4.6. Proof of 1.5.2.} Let $D=\prod_{i\in I} D_i$ be a product decomposition into isoclinic quasi-special $p$-divisible groups over $k$. For $i\in I$, let $c_i$ and $d_i$ be the codimension and the dimension (respectively) of $D_i$, and let $\alpha_i:={{d_i}\over {r_i}}$. Let $M=\bigoplus_{i\in I} M_i$ be the direct sum decomposition such that $(M_i,\phi)$ is the Dieudonn\'e module of $D_i$. As each $D_i$ is isoclinic, we have $n_{D_i}=\ell_{D_i}$ and $n_D=\ell_D$ (cf. Corollary 1.4.2). For $i\in I$ we have $\phi^{r_i}(M_i)=p^{d_i}M_i$, cf. Definition 1.5.1 (e). Let $m_i^\prime\in\dbN$ be the greatest divisor of $g.c.d.\{c_i,d_i\}$ such that for $(c_{i2},d_{i2},r_{i2}):=({{c_i}\over {m_i^\prime}},{{d_i}\over {m_i^\prime}},{{r_i}\over {m_i^\prime}})$ we have $\phi^{r_{i2}}(M_i)=p^{d_{i2}}M_i$. This identity implies that:

\medskip
{\bf (i)} we have $\alpha_{D_i}(r_{i2})=\beta_{D_i}(r_{i2})=d_{i2}$ and for all $q\in\dbN$ we have  $\alpha_{D_i}(q+r_{i2})=\alpha_{D_i}(q)+d_{i2}$ and $\beta_{D_i}(q+r_{i2})=\beta_{D_i}(q)+d_{i2}$.

\medskip
From (i) we get that for all $q\in\dbN$ we have $\delta_{D_i}(q+r_{i2})=\delta_{D_i}(q)$. From this and Proposition 4.3 (a) applied to $D_i$, we get that:

\medskip
{\bf (ii)} $n_{D_i}=\ell_{D_i}=\max\{\delta_{D_i}(q)|q\in\dbN\}=\max\{\delta_{D_i}(q)|q\in\{1,\ldots,r_{i2}\}\}$. 

\medskip
As the function $\beta_{D_i}(*)$ defined for $*\in\dbN$ is increasing, for all $q\in\{1,\ldots,r_{i2}\}$ we have $\delta_{D_i}(q)\le\beta_{D_i}(q)\le\beta_{D_i}(r_{i2})=d_{i2}$. From this and (ii) we get that $n_{D_i}=\ell_{D_i}\le d_{i2}$. It is easy to see that the $p$-divisible group $D_i^{\text{t}}$ is isoclinic quasi-special and that the analogue of the triple $(r_{i2},d_{i2},c_{i2})$ for it is $(r_{i2},c_{i2},d_{i2})$. Thus we have $n_{D_i^{\text{t}}}\le c_{i2}$. As $n_{D_i}=n_{D_i^{\text{t}}}$ (see Fact 4.2.1), we have $n_{D_i}\le c_{i2}$. Thus 
$$n_{D_i}=\ell_{D_i}\le\min\{c_{i2},d_{i2}\}\le\min\{c_i,d_i\}.\leqno (11a)$$ 
This proves Theorem 1.5.2 if $D=D_i$ i.e., if $I=\{i\}$.

We assume that $I$ has at least two elements. From Proposition 1.4.3 we get that $n_D=\ell_D\le\max\{1,n_{D_i}+n_{D_j}|i\in I,j\in I\setminus\{i\}\}$. From this and $(11a)$ we get that
$$n_D\le\max\{1,\min\{c_{i2}+c_{j2},d_{i2}+d_{j2}\}|i,j\in I,i\neq j\}.\leqno (11b)$$
As $c_{i2}+c_{j2}\le c$ and $d_{i2}+d_{j2}\le d$, we have $\min\{c_{i2}+c_{j2},d_{i2}+d_{j2}\}\le\min\{c,d\}$. From this and $(11b)$ we get that $n_D=\ell_D\le\max\{1,\min\{c,d\}\}$. But if $\min\{c,d\}=0$ (i.e., if $cd=0$), then $n_D=\ell_D=0$. Thus, regardless of what the product $cd$ is, we have $n_D=\ell_D\le\min\{c,d\}$. This ends the proof of Theorem 1.5.2.\endproof 

\medskip\noindent
{\bf 4.6.1. Scholium.} Let $\pi$ be a permutation of $J_r=\{1,\ldots,r\}$. Let $o$ be the order of $\pi$. We assume that $(D,\phi)$ is $(C_{\pi},\phi_{\pi})$; thus $D$ is $F$-cyclic and therefore (cf. Lemma 4.2.4 (a)) quasi-special. We will translate the property 4.6 (ii) and Proposition 4.4 (a) and (b) in terms only of the permutation $\pi$. Let $\pi=\prod_{i\in I} \pi_i$ be the product decomposition of the permutation $\pi$ into cycles. As in Subsection 1.5, we write $\pi_i=(e_{s_1},\ldots,e_{s_{r_i}})$ for some number $r_i\in\dbN$ (which can be $1$). Let $M_i$ be the $W(k)$-span of $\{e_{s_1},\ldots,e_{s_{r_i}}\}$. We have $pM_i\subseteq \phi_{\pi}(M_i)\subseteq M_i$, cf. the definition of $\phi_{\pi}$. Thus we have a direct sum decomposition $(M,\phi_{\pi})=\bigoplus_{i\in I} (M_i,\phi_{\pi})$ of Dieudonn\'e modules. Let $D=\prod_{i\in I} D_i$ be the product decomposition that corresponds to the direct sum decomposition $(M,\phi_{\pi})=\bigoplus_{i\in I} (M_i,\phi_{\pi})$. Each $p$-divisible group $D_i$ is $F$-circular and quasi-special. Let $c_i,d_i=r_i-c_i,\alpha_i={{d_i}\over {r_i}}\in\dbN\cup\{0\}$ be as in Subsection 1.5. 

Due to the property 4.6 (i), the difference $\delta_{D_i}(q)=\beta_{D_i}(q)-\alpha_{D_i}(q)$ depends only on $q$ modulo $o$. For $s\in J_r$ and $q\in\{1,\ldots,o\}$, let $\eta_q(s)\in\dbN\cup\{0\}$ be such that we have $\phi_{\pi}(e_s)=p^{\eta_q(s)}e_{\pi^q(s)}$. Thus $\eta_q(s)$ is the number of elements of the sequence $e_s,\pi(e_s),\ldots,\pi^{q-1}(e_s)$ that belong to the set $\{e_1,\ldots, e_d\}$. We have 
$$\alpha_{D_i}(q)=\min\{\eta_q(s_j)|j\in\{1,\ldots,r_i\}\}\;\; \text{and}\;\;\beta_{D_i}(q)=\max\{\eta_q(s_j)|j\in\{1,\ldots,r_i\}\}.\leqno (12a)$$ 
The $W(k)$-basis $\scrB=\{e_1,\ldots,e_r\}$ for $M$ is a disjoint union of $W(k)$-basis for $M_i$'s. We consider the standard $W(k)$-basis $\{e_s\otimes e_t^*|s,t\in J_r\}$ for $\End(M)$ defined by $\scrB$. We have $\phi_{\pi}^q(e_s\otimes e_t^*)=p^{\eta_q(s)-\eta_q(t)}e_{\pi^q(s)}\otimes e^*_{\pi^q(t)}$. If $\eta_o(s)>\eta_o(t)$ (resp. $\eta_o(s)=\eta_o(t)$ or $\eta_o(s)<\eta_o(t)$), then $e_s\otimes e_t^*$ belongs to $L_+$ (resp. to $L_0$ or $L_-$) and therefore the number $\ell(e_s,e_t)$ defined in Scholium 3.5.1 is $\max\{0,\eta_q(t)-\eta_q(s)|q\in\{1,\ldots,o\}\}$ (resp. is $\max\{0,\eta_q(s)-\eta_q(t)|q\in\{1,\ldots,o\}\}$). From Formula $(6a)$ we get that 
$$\ell_{C_{\pi}}=\max\{\eps_{C_{\pi}},\ell(e_s,e_t)|s,t\in J_r\}.\leqno (12b)$$
\noindent
{\bf 4.6.2. Example.} We assume that $c=d=8$; thus $r=16$. Let $\pi=\pi_1\pi_2$, where $\pi_1=(9\;10\;5\;11\;12\;6\;7\;8)$ and $\pi_2=(1\;2\;13\;3\;4\;14\;15\;16)$ are $8$-cycles. We have $o=8$. Let $D=C_{\pi}=D_1\oplus D_2$ be the product decomposition corresponding to the cycle decomposition $\pi=\pi_1\pi_2$. All Newton polygon slopes of $D_1$ and $D_2$ are ${1\over 2}$ and it is easy to see that $D_1^{\text{t}}$ is isomorphic to $D_1$. We have $(\delta_{D_1}(1),\ldots,\delta_{D_1}(8))=(1,2,2,2,2,2,1,0)$; thus $n_{D_1}=2$ (cf. property 4.6 (ii)). From Fact 4.2.1 we get that $n_{D_2}=2$. We have $(\eta_1(9),\ldots,\eta_8(9))=(0,0,1,1,1,2,3,4)$ and $(\eta_1(1),\ldots,\eta_8(1))=(1,2,2,3,4,4,4,4)$. Thus $(\eta_1(9)-\eta_1(1),\ldots,\eta_8(9)-\eta_8(1))=(-1,-2,-1,-2,-3,-2,-1,0)$. Thus $\ell(e_9,e_1)=\ell(e_1,e_9)=3$. This implies that $n_D=\ell_D\ge 3$. From Proposition 1.4.3 we get that $n_D\le 3$. Thus $n_{D_1\oplus D_2}=n_D=3=n_{D_1}+n_{D_2}-1$. 

Plenty of similar examples can be constructed in which the identity $n_{D_1\oplus D_2}=n_{D_1}+n_{D_2}-1$ holds and $D_1$ and $D_2$ are isoclinic of equal height and different dimension.

\medskip\smallskip\noindent
{\bf 4.7. Proof of 1.4.4.} The Dieudonn\'e module of $\tilde D$ is $(\tilde M,\phi)$, where $\tilde M$ is a $W(k)$-submodule of $M$ which contains $p^{\kappa}M$. Let $\tilde O=\tilde O_+\oplus\tilde O_0\oplus\tilde O_-$ be the level module of $(\tilde M,\phi)$. If $D$ and $\tilde D$ are ordinary, then the Proposition 1.4.4 is trivial. Thus to prove the Proposition 1.4.4, we can assume that $D$ and $\tilde D$ are not ordinary; thus from Remark 4.1.1 we get that $\ell_D$ (resp. $\ell_{\tilde D})$ is the smallest natural number such that we have $p^{\ell_D}\End(M)\subseteq O$ (resp. we have $p^{\ell_{\tilde D}}\End(\tilde M)\subseteq \tilde O$). As $p^{\kappa}M\subseteq \tilde M\subseteq M$, we have
$$p^{2\kappa}\End(M)\subseteq p^{\kappa}\End(\tilde M)\subseteq\End(M).\leqno (13)$$
For $q\in\dbN$ we have $\phi^q(p^{\kappa}\tilde O_+)\subseteq p^{\kappa}\tilde O_+\subseteq \End(\tilde M)\cap L_+$. As  $p^{\kappa}\tilde O_+\subseteq\End(M)$ (cf. (13)), we get that $p^{\kappa}\tilde O_+\subseteq O_+$. A similar argument shows that $p^{\kappa}\tilde O_0\subseteq O_0$ and $p^{\kappa}\tilde O_-\subseteq O_-$. Thus $p^{\kappa}\tilde O\subseteq O$. From this, the inclusion $p^{\ell_{\tilde D}}\End(\tilde M)\subseteq \tilde O$, and (13) we get that
$$p^{2\kappa+\ell_{\tilde D}}\End(M)\subseteq p^{\kappa+\ell_{\tilde D}}\End(\tilde M)\subseteq p^{\kappa}\tilde O\subseteq O\subseteq\End(M).$$
Thus $\ell_D\le 2\kappa+\ell_{\tilde D}$. Based on this inequality, the Proposition 1.4.4 follows from Corollary 1.4.2. This ends the proof of Proposition 1.4.4.\endproof 

\medskip\smallskip\noindent
{\bf 4.7.1. Example.} We assume that $c=d$. We have $r=2d$. Let $\pi:=(12\cdots r)$; its cyclic decomposition is $\pi=\pi_i$ (with $i$ as an index). As $\phi^r_{\pi}(M)=p^{r\over 2}(M)$, the $F$-circular $p$-divisible group $C_{\pi}$ is supersingular. If $d\ge 2$, then $\phi^2(M)\neq pM$ and therefore $C_{\pi}$ is not special. As $\phi_{\pi}^d(e_1)=p^de_{d+1}$ and $\phi_{\pi}^d(e_{d+1})=e_1$, we have $\alpha_D(d)=0$ and $\beta_D(d)=d$. This implies $\delta_D(d)=d$ and therefore from Proposition 4.3 (a) we get that $n_{C_{\pi}}=\ell_{C_{\pi}}\ge d$. As $n_{C_{\pi}}\le d$ (cf. Theorem 1.5.2), we have $n_{C_{\pi}}=d$. See [NV, Ex. 3.3] for a simpler proof that $n_{C_{\pi}}=d$ (in loc. cit. $C_{\pi}$ is denoted as $C_d$). Let $E$ be a supersingular $p$-divisible group over $k$ of height $2$. From [NV, Rm. 2.6 and Ex. 3.3] we get that the smallest number $\kappa\in\dbN\cup\{0\}$ such that we have an isogeny $C_{\pi}\twoheadrightarrow E^d$ is $\kappa:=\lceil{{d-1}\over 2}\rceil$. It is well known that $E^d$ is uniquely determined up to isomorphism by $E^d[p]$ (for instance, see [NV, Scholium 2.3] or see Formula $(12b)$ applied to the minimal permutation $(1\; d+1)\cdots (d\;r)$ of $J_r$). Thus $n_{E^d}=1$. If $d$ is odd, then $\kappa={{d-1}\over 2}$ and therefore $n_{C_{\pi}}=d=n_{E^d}+2\kappa$. This implies that in general, Proposition 1.4.4 is optimal.

\medskip\smallskip\noindent
{\bf 4.7.2. Example.} We assume that $r>0$ and that $D$ is isoclinic. Let $M_0$ and $\tilde M_0$ be as in Proposition 4.3 (c). Let $D_0$ and $\tilde D_0$ be the $p$-divisible groups over $k$ whose Dieudonn\'e modules are isomorphic to $(M_0,\phi)$ and $(\tilde M_0,\phi)$ (respectively), cf. Remark 4.3.1. To the inclusions $M_0\subseteq M$ and $\tilde M_0\subseteq M$ correspond isogenies $D\twoheadrightarrow D_0$ and $D\twoheadrightarrow\tilde D_0$ whose kernels are annihilated by $p^{u_D}=p^{\tilde u_D}$ and $p^{v_D}=p^{\tilde v_D}$ (respectively), cf. Proposition 4.3 (c). Let $j_D:=n_{D_0}$ and $\tilde j_D:=n_{\tilde D_0}$. From Propositions 1.4.4 and 4.3 (f) we get that:
$$u_D\le n_D=\ell_D\le \min\{j_D+2u_D,\tilde j_D+2v_D\}.\leqno (14a)$$
If $(c_1,d_1,r_1)$ is as in Definition 4.1 (b), then $j_D\le \min\{c_1,d_1\}$ (cf. $(11a)$). From this and $(14a)$ we get:
$$u_D\le n_D\le 2u_D+\min\{c_1,d_1\}.\leqno (14b)$$

\bigskip
\noindent
{\boldsectionfont 5. On the Main Theorem B}
\bigskip

In Subsection 5.1 we prove the Main Theorem B. Subsections 5.2 and 5.3 present two applications of the Main Theorem B. For instance, Theorem 5.3 presents applications to extensions between two minimal $p$-divisible groups over $k$. We recall that $(M,\phi)$ is the Dieudonn\'e module of $D$.

\medskip\smallskip\noindent
{\bf 5.1. The proof of the Main Theorem B.} If $n_D\le 1$, then $D[p]$ is minimal (cf. Definition 1.5.1 (d)). If $D[p]$ is minimal, then there exists a $p$-divisible group $\tilde D$ over $k$ such that $n_{\tilde D}\le 1$ and $\tilde D[p]$ is isomorphic to $D[p]$; the codimension and the dimension of $\tilde D$ are $c$ and $d$ (respectively) and thus from the very definition of $n_{\tilde D}$ we get that $D$ is isomorphic to $\tilde D$ and therefore that we have $n_D=n_{\tilde D}\le 1$. Thus we have $n_D\le 1$ if and only if $D[p]$ is minimal. As $n_D\Le \ell_D$ (see Corollary 1.4.2),  1.6 (a) implies 1.6 (b). Thus to end the proof of the Main Theorem B, it suffices to show that 1.6 (b) implies 1.6 (c) and that 1.6 (c) implies 1.6 (a). 

\medskip\noindent
{\bf 5.1.1. On $1.6\, (b)\Rightarrow 1.6\, (c)$.} Let $\phi_1,\vartheta_1:M/pM\to M/pM$ be the reductions modulo $p$ of $\phi,\vartheta:M\to M$. In [Kr] (see also [Oo1, Subsect. (2.3) and Lem. (2.4)] and [Mo, Subsect. 2.1]) it is shown that there exists a $k$-basis $\{b_1,\ldots,b_r\}$ for $M/pM$ and a permutation $\pi$ of $J_r=\{1,\ldots,r\}$ such that the following two properties hold:
\medskip
{\bf  (i)} if $s\in\{1,\ldots,d\}$, then $\phi_1(b_s)=0$ and $\vartheta_1(b_{\pi(s)})=b_s$, and 
\smallskip
{\bf (ii)} if $s\in\{d+1,\ldots,r\}$, then $\phi_1(b_s)=b_{\pi(s)}$ and $\vartheta_1(b_{\pi(s)})=0$.
\medskip
Let $\vartheta_{\pi}:=p\phi_{\pi}^{-1}:M\to M$; if $s\in\{1,\ldots,d\}$, then $\vartheta_{\pi}(e_{\pi(s)})=e_s$, and if $s\in\{d+1,\ldots,r\}$, then $\vartheta_{\pi}(e_{\pi(s)})=pe_s$. Properties (i) and (ii) imply that the $k$-linear map $M/pM\to M/pM$ that takes $b_s$ to $e_s$ modulo $p$, is an isomorphism between  $(M/pM,\phi_1,\vartheta_1)$ and the reduction modulo $p$ of $(M,\phi_{\pi},\vartheta_{\pi})$. This means that $D[p]$ is isomorphic to $C_{\pi}[p]$, cf. the classical Dieudonn\'e theory. As $n_D\le 1$, we get that $D$ is isomorphic to $C_{\pi}$. 

We check that $\pi$ is a minimal permutation in the sense of Definition 1.5.1 (b). Let $\pi=\prod_{i\in I} \pi_i$ be the product decomposition of $\pi$ into cycles. We write $\pi_i=(e_{s_1},\ldots,e_{s_{r_i}})$, where $r_i\in\dbN$. Let $c_i$, $d_i=r_i-c_i$, and $\alpha_i={{d_i}\over {r_i}}$ be as in Subsection 1.5. Let $M_i:=\bigoplus_{u=1}^{r_i} W(k)e_{s_u}$. Let $D=\prod_{i\in I} D_i$ be the product decomposition defined by the direct sum decomposition $(M,\phi_{\pi})=\bigoplus_{i\in I} (M_i,\phi_{\pi})$. Each $D_i$ is an $F$-circular $p$-divisible group over $k$ and therefore isoclinic. From Proposition 4.4 (b) we get that $n_{D_i}\le n_D\le 1$. But $n_{D_i}=\ell_{D_i}=\max\{\delta_{D_i}(q)|q\in\dbN\}$, cf. Corollary 1.4.2 and Proposition 4.3 (a). From the last two sentences we get that for all $q\in\dbN$ we have $\delta_{D_i}(q)\in\{0,1\}$. Thus either $\alpha_{D_i}(q)=\beta_{D_i}(q)$ or $\alpha_{D_i}(q)+1=\beta_{D_i}(q)$. If $\alpha_{D_i}(q)=\beta_{D_i}(q)$, then from Lemma 4.2.3 (a) we get that $\alpha_{D_i}(q)=\beta_{D_i}(q)=q\alpha_i$. If $\alpha_{D_i}(q)+1=\beta_{D_i}(q)$, then from Lemma 4.2.3 (a) we get that either $(\alpha_{D_i}(q),\beta_{D_i}(q))=([q\alpha_i],[q\alpha_i]+1)$ or $(\alpha_{D_i}(q),\beta_{D_i}(q))=(q\alpha_i-1,q\alpha_i)$. But the second possibility is excluded by Lemma 4.2.3 (b). We conclude that in all cases we have $\alpha_{D_i}(q),\beta_{D_i}(q)\in \{[q\alpha_i],[q\alpha_i]+1\}$. Therefore $p^{[q\alpha_i]+1}M_i\subseteq\phi^q(M_i)\subseteq p^{[q\alpha_i]}M_i$. Thus for each $u\in\{1,\ldots,r_i\}$, we have $\phi_{\pi}^q(e_{s_u})=p^{[q\alpha_i]+\eps_q(s_u)}e_{\pi^q(s_u)}$ for some number $\eps_q(s_u)\in\{0,1\}$. As this property holds for all pairs $(q,i)\in\dbN\times I$, $\pi$ is a minimal permutation. As $D$ is isomorphic to $C_{\pi}$, we get that $D$ is minimal. Thus 1.6 (b) implies 1.6 (c).

\medskip\noindent
{\bf 5.1.2. On $1.6\, (c)\Rightarrow 1.6\, (a)$.} To prove that 1.6 (c) implies 1.6 (a), we can assume that $\pi$ is a minimal permutation of $J_{r}$, that $D=C_{\pi}$,  and that $\phi=\phi_{\pi}$. Let $\pi=\prod_{i\in I} \pi_i$, $M=\bigoplus_{i\in I} M_i$, and $D=\prod_{i\in I} D_i$ be the decompositions obtained as in Subsubsection 5.1.1. For $i\in I$, let $\pi_i=(e_{s_1},\ldots,e_{s_{r_i}})$, $c_i$, $d_i=r_i-c_i$, and $\alpha_i={{d_i}\over {r_i}}$ be as in Subsection 1.5. As the permutation $\pi$ is minimal, for all $u\in\{1,\ldots,r_i\}$ we have $\phi_{\pi}^q(e_{s_u})=p^{[q\alpha_i]+\eps_q(s_u)}e_{\pi^q(s_u)}$ for some number $\eps_{q}(s_u)\in\{0,1\}$. This implies that $p^{[q\alpha_i]+1}M_i\subseteq \phi_{\pi}^q(M_i)\subseteq p^{[q\alpha_i]}M_i$. Thus 
$$\alpha_{D_i}(q),\beta_{D_i}(q)\in\{[q\alpha_i],[q\alpha_i]+1\}.\leqno (15)$$ 
From (15) and the fact that $\delta_{D_i}(q)\ge 0$, we get that $\delta_{D_i}(q)\in\{0,1\}$. From this and Proposition 4.3 (a), we get that $n_{D_i}=\ell_{D_i}\le 1$. If $D=D_i$ (i.e., if $I=\{i\}$), then $\ell_D\le 1$ and thus 1.6 (a) holds. If $I$ has at least two elements, then from Proposition 1.4.3 we get that $\ell_D=n_D\le 1$. Thus regardless of what $I$ is, we have $\ell_D\le 1$. This ends the argument that the implication $1.6\, (c)\Rightarrow  1.6\, (a)$ holds. This ends the proof of the Main Theorem B.\endproof 

\bigskip\noindent
{\bf 5.2. Corollary.} {\it We assume that $\ell_D\le 2$. Then $n_D=\ell_D$.}

\medskip
\proof
If $n_D\le 1$, then $D$ is minimal (cf. Main Theorem B) and therefore $F$-cyclic; thus $n_D=\ell_D$ (cf. Theorem 1.5.2). As $n_D\le\ell_D\le 2$, we have $n_D=\ell_D$ even if $n_D=2$.\endproof 

\medskip
The next Theorem generalizes and refines [Va1, Prop. 4.5.1]. 

\medskip\smallskip\noindent
{\bf 5.3. Theorem.} {\it We assume that we have a short exact sequence $0\to D_1\to D\to D_2\to 0$ of $p$-divisible groups over $k$, with $D_1$ and $D_2$ as minimal $p$-divisible groups. 

\medskip
{\bf (a)} Then we have $n_D\le\ell_D\le 3$.

\smallskip
{\bf (b)} We assume that $d=c\ge 3$ and that $D_1$ and $D_2$ are isoclinic of Newton polygon slopes ${1\over d}$ and ${{d-1}\over d}$ (respectively). Then $n_D=\ell_D\le2$.}

\medskip
\proof
Let $0\to D_1\to \tilde D\to D_2\to 0$ be the pull forward of our initial short exact sequence via the multiplication by $p$ isogeny $D_1\twoheadrightarrow D_1$. The kernel of the resulting isogeny $D\twoheadrightarrow\tilde D$ is annihilated by $p$. The $p$-divisible group $D_{12}:=D_1\oplus D_2$ is minimal. As $\tilde D[p]$ is isomorphic to $D_{12}[p]$, from the equivalence between 1.6 (b) and (c) we get that $\tilde D$ is isomorphic to $D_{12}$ and that $n_{\tilde D}\le 1$. As $\tilde D$ is $F$-cyclic and thus a direct sum of isoclinic $p$-divisible groups, we have $n_D\le\ell_D\le n_{\tilde D}+2\le 3$ (cf. Proposition 1.4.4). Thus (a) holds.

We prove (b). We know that there exists an isogeny $D_{12}\twoheadrightarrow D$ whose kernel $K$ is annihilated by $p$. We will choose such an isogeny of the smallest degree possible. It is well known that up to isomorphisms, there exists a unique $p$-divisible group over $k$ of height $d$ and Newton polygon slope ${*\over d}$, where $*\in\{1,d-1\}$ (see [De, Ch. IV, Sect. 8]). Thus if $K$ has a proper subgroup scheme $K_1$ (resp. $K_2$) whose image in $D_1$ (resp. in $D_2$) is trivial, then $D_{12}^\prime:=D_{12}/K_1$ (resp. $D_{12}^\prime:=D_{12}/K_2$) is isomorphic to $D_{12}$ and thus we would get an isogeny $D_{12}\arrowsim D_{12}^\prime\twoheadrightarrow D$ of smaller degree. This implies that the projections of $K$ on $D_1$ and $D_2$ are monomorphisms. Thus the codimension and the dimension of $K$ are both at most $1$. Based on the last two sentences, as $d\ge 3$ we easily get that $K$ is either trivial or isomorphic to $\pmb{\alpha}_p$. If $K$ is trivial, then $D$ is minimal and therefore we have $n_D\le 1$ (in fact we have $n_D=1$). Thus to prove (b), we can assume that $K$ is isomorphic to $\pmb{\alpha}_p$. We reached the case when we have isogenies 
$$D_{12}\twoheadrightarrow D_{12}/\pmb{\alpha}_p\arrowsim D\twoheadrightarrow D_{12}/(\pmb{\alpha}_p\times_k\pmb{\alpha}_p)=D_1/\pmb{\alpha}_p\times_k D_2/\pmb{\alpha}_p.$$ 
At the level of Dieudonn\'e modules, this means the following things. Let $N_{12}:=\bigoplus_{s=1}^{2d} W(k)e_s$ be a free $W(k)$-module of rank $r=2d$. Let $\phi:N_{12}\to N_{12}$ be the $\sigma$-linear endomorphism such that it takes $(e_1,\ldots,e_d)$ and $(e_{d+1},\ldots,e_{2d})$ to $(pe_2,e_3,\ldots,e_d,e_1)$ and $(pe_{d+2},\ldots,pe_{2d},e_{d+1})$ (respectively). We can assume that $(N_{12},\phi)$ is the Dieudonn\'e module of $D_1/\pmb{\alpha}_p\times_k D_2/\pmb{\alpha}_p$ (cf. the mentioned uniqueness property). As $D\twoheadrightarrow D_{12}/(\pmb{\alpha}_p\times_k\pmb{\alpha}_p)=D_1/\pmb{\alpha}_p\times_k D_2/\pmb{\alpha}_p$ is an isogeny of kernel $\pmb{\alpha}_p$ and as $K$ maps monomorphically to both $D_1$ and $D_2$, there exists an invertible element $\gamma\in W(k)$ such that we can identify $M$ with $N_{12}+W(k)({{\gamma}\over p}e_1+{1\over p}e_{d+1})$. Moreover, if $M_{12}:=N_{12}+W(k){1\over p}e_1+W(k){1\over p}e_{d+1}$, then $(M_{12},\phi)$ is the Dieudonn\'e module of $D_{12}$. 

We check that $\ell_D\le 2$. We have 
$$p^2\End(M)\subseteq p\Hom(M_{12},N_{12})+W(k)[(p\gamma e_1+pe_{d+1})\otimes e_1^*]\subseteq\End(M).\leqno (16)$$ 
The latticed $F$-isocrystal $(\Hom(M_{12},N_{12}),\phi)$ is isomorphic to $\End(M_{12},\phi)$ and moreover we have $\ell_{D_{12}}\le 1$. From this and  the fact that $\Hom(M_{12},N_{12})\subseteq\End(M)$, we get that $O$ contains $p\Hom(M_{12},N_{12})$. It is easy to see that for all $q\in\dbN$ we have $\phi^q((p\gamma e_1+pe_{d+1})\otimes e_1^*)\in\End(M)$; like $\phi((p\gamma e_1+pe_{d+1})\otimes e_1^*)=(p\sigma(\gamma) e_2+pe_{d+2})\otimes e_2^*$, $\phi^2((p\gamma e_3+pe_{d+1})\otimes e_1^*)=(p\sigma^2(\gamma) e_3+p^2e_{d+3})\otimes e_3^*$, etc. Thus $(p\gamma e_1+pe_{d+1})\otimes e_1^*\in O_+\oplus O_0$, cf. Lemma 2.4.  Based on (16) we conclude that $p^2\End(M)\subseteq O$. Thus $\ell_D\le 2$. From Corollary 5.2 we get that $n_D=\ell_D\le 2$. Thus (b) holds.\endproof

\bigskip\noindent
{\bf Acknowledgments.} We would like to thank University of Arizona and Binghamton University for providing good conditions during the writing of this paper. We would like to thank M.-H. Nicole for introducing us to Traverso's work on February 2006, for many valuable comments, and for suggesting the ideas of the second paragraph of the proof of Theorem 5.3 (cf. [Ni, Thm. 1.5.2]). We would like to thank the referee for many valuable comments and suggestions. 

\bigskip
\references{37}
{\nspace{

\Ref[De]
M. Demazure, 
\sl Lectures on $p$-divisible groups, 
\rm Lecture Notes in Math., Vol. {\bf 302}, Springer-Verlag, Berlin-New York, 1972.

\Ref[Di] 
J. Dieudonn\'e, 
\sl Groupes de Lie et hyperalg\`ebres de Lie sur un corps de caract\'erisque $p>0$. VII, 
\rm Math. Annalen {\bf 134} (1957), pp. 114--133.

\Ref[Kr] H. Kraft, 
\sl Kommutative algebraische p-Gruppen (mit Anwendungen auf p-divisible Gruppen und abelsche Variet\"aten), 
\rm manuscript 86 pages, Univ. Bonn, 1975.

\Ref[Ma]
J. I. Manin, 
\sl The theory of formal commutative groups in finite characteristic, 
\rm Russian Math. Surv. {\bf 18} (1963), no. 6, pp. 1--83.

\Ref[Mo]
B. Moonen,
\sl Group schemes with additional structures and Weyl group cosets,
\rm Moduli of abelian varieties (Texel Island, 1999), pp. 255--298, Progr. Math., Vol. {\bf 195}, Birkh\"auser, Basel, 2001.

\Ref[Ni]
M.-H. Nicole,
\sl Superspecial abelian varieties, theta
series and the Jacquet-Langlands correspondence, 
\rm Ph.D. thesis, McGill University,  October 2005.

\Ref[NV]
M.-H. Nicole and A. Vasiu,
\sl Minimal truncations of supersingular $p$-divisible groups,
\rm Indiana Univ. Math. J. {\bf 56} (2007), no. 6, pp. 2887--2897.

\Ref[Oo1] 
F. Oort, 
\sl A stratification of a moduli space of abelian varieties, 
\rm Moduli of abelian varieties (Texel Island, 1999), pp. 345--416, Progr. Math., Vol. {\bf 195}, Birkh\"auser, Basel, 2001.

\Ref[Oo2] 
F. Oort, 
\sl Foliations in moduli spaces of abelian varieties, 
\rm J. Amer. Math. Soc. {\bf 17} (2004), no. 2, pp. 267--296.

\Ref[Oo3]
F. Oort,
\sl  Minimal $p$-divisible groups.  
\rm Ann. of Math. (2)  {\bf 161}  (2005),  no. 2, pp. 1021--1036. 

\Ref[Oo4]
F. Oort,
\sl Simple $p$-kernels of $p$-divisible groups,
\rm  Adv. Math. {\bf 198}  (2005),  no. 1, pp. 275--310.

\Ref[Tr1]
C. Traverso,
\sl Sulla classificazione dei gruppi analitici commutativi di caratteristica positiva,
\rm Ann. Scuola Norm. Sup. Pisa (3) {\bf 23} (1969), no. 3, pp. 481--507.

\Ref[Tr2]
C. Traverso,
\sl p-divisible groups over fields, 
\rm Symposia Mathematica, Vol. {\bf XI} (Convegno di Algebra Commutativa, INDAM,
Rome, 1971), pp. 45--65, Academic Press, London, 1973.

\Ref[Tr3]
C. Traverso,
\sl Specializations of Barsotti--Tate groups,
\rm  Symposia Mathematica, Vol. {\bf XXIV} (Sympos., INDAM, Rome, 1979), pp. 1--21, Acad. Press, London-New York, 1981.  

\Ref[Va1]
A. Vasiu,
\sl Crystalline Boundedness Principle,
\rm Ann. Sci. \'Ecole Norm. Sup. {\bf 39} (2006), no. 2, pp. 245--300.

\Ref[Va2]
A. Vasiu,
\sl Generalized Serre--Tate ordinary theory,
\rm http://arxiv.org/abs/math/0208216.

\Ref[Va3]
A. Vasiu,
\sl Mod p classification of Shimura F-crystals,
\rm http://arxiv.org/abs/math/0304030.

}}
\noindent
\hbox{Adrian Vasiu,}
\hbox{Department of Mathematical Sciences, Binghamton University,}
\hbox{Binghamton, New York 13902-6000, U.S.A.}
\hbox{e-mail: adrian\@math.binghamton.edu}

\enddocument